\documentclass{amsart}

\usepackage{amsmath, amssymb}
\usepackage{amsfonts}
\usepackage[all]{xy}

\usepackage{epic}
\usepackage{graphics}
\usepackage{epsfig}

\usepackage{picture}

\usepackage{epstopdf}

\numberwithin{equation}{section}

\DeclareMathOperator{\Aut}{Aut}

\DeclareMathOperator{\id}{id}
\DeclareMathOperator{\Hom}{Hom}
\DeclareMathOperator{\dhom}{hom}
\DeclareMathOperator{\Br}{Br}

\DeclareMathOperator{\GL}{GL}
\DeclareMathOperator{\rank}{rank}

\DeclareMathOperator{\Sph}{Br} 

\DeclareMathOperator{\GM}{GM}
\DeclareMathOperator{\Sim}{Sim}
\DeclareMathOperator{\Mod}{mod-}
\DeclareMathOperator{\Per}{Per}
\DeclareMathOperator{\Cone}{Cone}
\DeclareMathOperator{\Int}{Int}

\def\bs{\backslash}
\def\lam{\lambda}
\def\endproof{\hfill{$\square$} \medskip}

\def\C{\mathbb{C}}
\def\R{\mathbb{R}}

\def\Z{\mathbb{Z}}
\def\alp{\alpha}
\def\sli{\mathcal{P}}
\def\H{\mathcal{H}}
\def\class{\mathcal{C}}
\def\Stab{\mathrm{Stab}}

\def\D{\mathcal{D}}
\def\P{\mathbb{P}}
\def\Q{Q(N,n)}
\def\Crit{\mathrm{Crit}}
\def\Zero{\mathrm{Zero}}
\def\lto{\longrightarrow}
\def\GQ{\Gamma_N Q}
\def\DQ{\D^N_Q}
\def\iso{\xrightarrow{\sim}}

\def\xto{\xrightarrow}
\def\Ang{\mathrm{A}}
\def\EG{\mathrm{EG}}
\def\QA{\overrightarrow{A_n}}

\def\d{d_n^N}
\def\Ho{\textbf{H}}
\def\LQ{L_{\mathcal{Q}}}
\def\CQ{\mathcal{Q}}
\def\CY{\mathrm{CY}}

\newtheorem{thm}{Theorem}[section]
\newtheorem{prop}[thm]{Proposition}
\newtheorem{cor}[thm]{Corollary}
\newtheorem{rem}[thm]{Remark}
\newtheorem{lem}[thm]{Lemma}
\newtheorem{conjecture}[thm]{Conjecture}

\theoremstyle{definition}
\newtheorem{defi}[thm]{Definition}

\renewcommand{\Im}{\text{Im\,}}
\renewcommand{\Re}{\text{Re\,}}
\newcommand{\norm}[1]{\lVert #1 \rVert}
\newcommand{\abv}[1]{\lvert #1 \rvert}

\author{Akishi Ikeda}
\address{Kavli Institute for the Physics and Mathematics of the Universe (WPI), UTIAS, 
The University of Tokyo, Kashiwa, Chiba 277-8583, Japan}
\email{akishi.ikeda@ipmu.jp}

\begin{document}

\title[Stability conditions on $\text{CY}_N$ categories associated 
to $A_n$-quivers]
{Stability conditions on $\text{CY}_N$ categories associated to 
$A_n$-quivers and period maps}

\begin{abstract}

In this paper, we study the space of stability conditions 
on a certain $N$-Calabi-Yau ($\CY_N$) category associated to an  $A_n$-quiver. 
Recently, 
Bridgeland and Smith constructed 
stability conditions on some $\CY_3$ categories from meromorphic quadratic differentials with simple zeros. 
Generalizing their results to higher dimensional 
Calabi-Yau categories, we describe the space of stability 
conditions as the universal 
cover of the space of polynomials of degree $n+1$ 
with simple zeros.    
In particular, central charges of stability conditions on $\CY_N$ 
categories are constructed  
as the periods of quadratic differentials with zeros of order $N-2$ 
which are associated to polynomials. 
\end{abstract}

\maketitle

\section{Introduction}
The notion of a stability condition on a triangulated category 
$\D$ was 
introduced by Bridgeland \cite{Br1} following Douglas's work 
on $\Pi$-stability \cite{Dou}. Bridgeland also showed that 
the space of stability conditions $\Stab(\D)$ has 
the structure of a complex manifold. 

In this paper, we study the space of stability conditions 
on a certain $N$-Calabi-Yau ($\CY_N$) category $\D_n^N$ 
associated to an $A_n$-quiver \cite{ST,T}. 
The category 
$\D_n^N$ is defined by  
the derived category of finite dimensional dg modules 
over the Ginzburg $\CY_N$ dg algebra \cite{Ginz} 
of the $A_n$-quiver. 

We show that 
the distinguished connected component 
$\Stab^{\circ}(\D_n^N) \subset \Stab(\D_n^N)$ 
is isomorphic to the universal cover 
of the space of polynomials  
$p_n(z)=z^{n+1}+u_1 z^{n-1}+\cdots + u_n\,(u_i\in \C)$ 
with simple zeros. 
The central charges of stability conditions in $\Stab^{\circ}(\D_n^N)$ 
are constructed as the periods of quadratic differentials of 
the form 
$p_n(z)^{N-2} dz^{\otimes 2}$ on the Riemann sphere $\P^1$. 

To prove the result, we use the recent results 
of Bridgeland and Smith  \cite{BrSm}. Motivated by the work of 
physicists Gaiotto-Moore-Neitzke \cite{GMN}, 
they constructed stability conditions 
on a class of $\CY_3$ categories associated to triangulations of 
marked bordered surfaces
from quadratic differentials with simple zeros. 

By considering quadratic differentials with zeros of order $N-2$ 
instead of simple zeros, 
we can generalize the framework of \cite{BrSm} 
for $\CY_3$ categories to $\CY_N$ categories.  
Here we only treat the simplest case that 
the surface is a disk with some 
marked points on the boundary. In the future, we expect that 
the framework of \cite{BrSm} can be  
generalized to a class of $\CY_N$ categories 
associated to $N$-angulations of marked bordered surfaces, 
though these categories are not (well) defined yet.

\subsection{Background geometry}
Let $N \ge 2$ be an integer. 
Here we review the geometric background for this paper and 
explain why we consider 
the periods of the differential form $p_n(z)^{\frac{N-2}{2}}dz$. 
Section 2 in \cite{T} is 
an excellent survey of related topics. 

The category $\D_n^N$ was
first introduced as the derived Fukaya category \cite{KhSe,T} 
of Lagrangian $N$-spheres which form the basis of vanishing cycles 
in the universal deformation of the $A_n$-singularity
\begin{equation*}
X^N_p:=
\{x_1^2 +x_2^2+\dots +x_N^2=p_n(z)\} \subset \C^N \times \C
\end{equation*}
where $p_n(z)=z^{n+1}+u_1 z^{n-1}+ \cdots +u_n\,
(u_i \in \C)$ 
with simple zeros. We consider a 
symplectic structure on $X^N_p$ coming    
from the standard K\"ahler form on $\C^N \times \C$, 
and a nowhere vanishing holomorphic $N$-form 
$\Omega_p$ on $X^N_p$ given by taking Poincar\'e residue 
of the standard form $dx_1 \cdots dx_N dz$ on $\C^{N+1}$.  

Consider the projection 
\begin{equation*}
\pi \colon X^N_p \to \C,\quad (x_1,\dots,x_N,z) \mapsto z.
\end{equation*}
If $z \in \C$ is not a root of $p_n$, then the fiber $\pi^{-1}(z)$ 
is an affine quadric of complex dimension $N-1$ with a real slice $S^{N-1} \subset \pi^{-1}(z)$.  
Thus, each Lagrangian $N$-sphere is constructed 
by the total space of the $S^{N-1}$-fibration over a path 
$\gamma \colon [0,1] \to \C$ which connects 
distinct roots of $p_n$ 
(see \cite[Section 6]{TY} or \cite[Section 2]{T}). 

As in \cite{Dou}, a central charge ``in the physics literature'' of 
a Lagrangian $L \subset X^N_p$ is given by the period 
\begin{equation*}
\int_L \Omega_p. 
\end{equation*}
Assume that $L$ is a $S^{N-1}$-fibration over 
the path $\gamma \colon [0,1] \to \C$. 
Then, up to constant multiplication, the above integration is reduced 
to the integration of the differential form 
$p_n(z)^{\frac{N-2}{2}}dz$ along the path $\gamma$ 
(see \cite[Section 6]{TY}):
\begin{equation*}
\int_{\gamma}p_n(z)^{\frac{N-2}{2}}dz.
\end{equation*}
From this observation, we expect that the central charges of 
Bridgeland stability conditions on $\D_n^N$ 
can be constructed by using 
the period of the differential form $p_n(z)^{\frac{N-2}{2}}dz$. 
 
\subsection{Summary of results}
First we explain the space of polynomials and 
the associated period map (Section \ref{sec_2}). 
Let $M_n$ be the space of polynomials $p_n(z)=z^{n+1}+u_1
z^{n-1}+\cdots+ u_n\,(u_i \in \C)$ 
with $n+1$ simple zeros. 
Note that the fundamental group of $M_n$ is the 
Artin group $B_{n+1}$. The action of $\C^*$ on $M_n$ 
is given by $(u_1,\dots,u_n) \mapsto (ku_1,\dots,k^n u_n)$ 
for $k \in \C^*$.

We regard $p_n(z)$ as a meromorphic function on  
the Riemann sphere $\P^1 =\C \cup \{\infty\}$. 
For $p_n(z) \in M_n$, we define the quadratic differential 
on $\P^1$ by $\phi(z):=p_n(z)^{N-2} dz^{\otimes 2}$ and 
denote by $\Q$ the space of such differentials. 

To consider the period of $\sqrt{\phi}=p_n(z)^{\frac{N-2}{2}}dz$, 
we need to introduce homology groups $H_+(\phi)$ and 
$H_-(\phi)$ which provide appropriate cycles for the 
integration of $\sqrt{\phi}$.
If $N$ is even, 
the homology group is defined by the relative homology group 
$H_+(\phi)=H_1(\C,\Zero(\phi);\Z)$ where $\Zero(\phi)\subset \C$ 
is the set of zeros of $\phi$. If $N$ is odd, the homology 
group is defined by $H_-(\phi):=H_1(S \bs \pi^{-1}(\infty);\Z)$ 
where $\pi \colon S \to \P^1$ is the double cover of 
the hyperelliptic curve $S=\{y^2=p_n(z)\}$. 

A linear map $Z_{\phi} \colon H_{\pm}(\phi) \lto \C$, called 
the period of $\phi$, is defined by
\begin{equation*}
Z_{\phi}(\gamma):=\int_{\gamma}\sqrt{\phi}
\end{equation*}    
for $\gamma \in H_{\pm}(\phi)$. 

Let $\Gamma$ be a free abelian group of rank $n$. 
An isomorphism of abelian groups 
$\theta \colon \Gamma \xto{\sim} H_{\pm}(\phi)$ is called 
a $\Gamma$-framing of $\phi \in \Q$. 
Denote by $\Q^{\Gamma}$ the set of framed differentials. 
For a framed differential $(\phi,\theta) \in \Q^{\Gamma}$, 
the composition of  $\theta \colon \Gamma \to H_{\pm}(\phi)$ 
and $Z_{\phi} \colon H_{\pm}(\phi) \to \C$ gives a 
linear map $Z_{\phi} \circ \theta \colon \Gamma \to \C$. 
We define the period map by
\begin{equation*}
\mathcal{W}_N \colon \Q^{\Gamma} 
\to \Hom_{\Z}(\Gamma ,\C),\quad 
(\phi,\theta) \mapsto Z_{\phi} \circ \theta.
\end{equation*} 

Fix a point $*=(\phi_0,\theta_0) \in \Q^{\Gamma}$ 
and let $\Q^{\Gamma}_* \subset \Q^{\Gamma}$ 
be the connected component 
containing  $*$. Then the universal 
cover of $\Q^{\Gamma}_* $ is 
isomorphic to the universal cover $\widetilde{M_n}$ of $M_n$. 
Thus the period map $\mathcal{W}_N$ is extended on 
$\widetilde{M_n}$:
\begin{equation*}
\mathcal{W}_N \colon \widetilde{M_n} \lto \Hom(\Gamma,\C).  
\end{equation*}

Next we consider the category 
$\D_n^N$ (Section \ref{sec_4}) 
and stability conditions on it (Section \ref{sec:stab}). 
For a quiver $Q$, we can define 
a dg algebra $\Gamma_N Q$, 
called the Ginzburg $\CY_N$ dg algebra \cite{Ginz}. Let $\D_{\mathrm{fd}}(\Gamma_N Q)$ be the derived category 
of dg modules over $\Gamma_N Q$ with finite total dimension. 
It was proved by Keller \cite{Kel2} that the category 
$\D_{\mathrm{fd}}(\Gamma_N Q)$ is a $\CY_N$ 
triangulated category. Let $Q=\QA$ be the $A_n$-quiver 
and set $\D^N_n:=\D_{\mathrm{fd}}(\Gamma_N \QA)$.
An object $S \in \D^N_n$ is called $N$-spherical if 
\begin{align*}
\Hom_{\D^N_n}(S,S[i]) = 
\begin{cases}
k \quad \text{if} \quad i=0,N  \\
\,0 \quad \text{otherwise}. 
\end{cases}
\end{align*}
There is a full subcategory $\H_0$, called 
the standard heart of $\D_n^N$, and there are 
$N$-spherical objects $S_1,\dots,S_n \in \H_0$ such that 
$\H_0$ is generated by these $N$-spherical objects. 

From such objects, Seidel and Thomas defined autoequivalences 
$\Phi_{S_1},\dots,\Phi_{S_n} \in \Aut(\D_n^N)$, called 
spherical twists \cite{ST}. They also proved that 
the subgroup $\Sph(\D_n^N) \subset \Aut(\D_n^N)$ generated 
by $\Phi_{S_1},\dots,\Phi_{S_n}$ is isomorphic to the braid group 
$B_{n+1}$.

A stability condition $\sigma=(Z,\H)$ on a triangulated 
category $\D$ consists of a heart $\H \subset \D$ 
of a bounded t-structure 
and a linear map $Z \colon K(\H) \to \C$ called 
the central charge, which satisfy some axioms. 
In \cite{Br1}, Bridgeland showed 
that the 
set of stability conditions $\Stab(\D)$ has the structure of 
a complex manifold, and the projection map 
\begin{equation*}
\mathcal{Z} \colon \Stab(\D) \lto \Hom_{\Z}(K(\D),\C),\quad 
(Z,\H) \mapsto Z
\end{equation*}
is a local isomorphism. 
There are two group actions on $\Stab(\D)$; 
one is the action of $\C$ and the other 
is the action of $\Aut(\D)$.  

Consider the category $\D^N_n$. There is a distinguished connected 
component $\Stab^{\circ}(\D_n^N) \subset \Stab(\D_n^N)$ which 
contains stability conditions with the standard 
heart $\H_0$. The action of the 
group $\Sph(\D^N_n)$ preserves 
the connected component $\Stab^{\circ}(\D_n^N) $.  
Since $\Sph(\D_n^N) \cong B_{n+1}$, the space 
$\Stab^{\circ}(\D_n^N)$ has the action of 
the braid group $B_{n+1}$. 

Our main result is the following. 

\begin{thm}[Theorem \ref{thm_main}]
Assume that $N \ge 2$. 
There is a $B_{n+1}$-equivariant and 
$\C$-equivariant isomorphism of complex manifolds 
$\widetilde{K_N}$ such that the diagram
\begin{equation*}
\xymatrix{
\widetilde{M_n} \ar[rr]^{\widetilde{K_N}} \ar[dr]_{\mathcal{W}_N} && 
\Stab^{\circ}(\D^N_n)  \ar[dl]^{\mathcal{Z}}  \\
& \Hom_{\Z}(\Gamma,\C) &
}
\end{equation*}
commutes. 
\end{thm}
When $N=2$, the result was proved by Thomas \cite{T}. 
When $N=3$ and $n=2$, the result was proved by 
Sutherland \cite{Suther}. The result for  $N=3$ and arbitrary $n$ 
follows from the result of 
Bridgeland-Smith \cite{BrSm} with some modifications 
(see Example 12.1 in \cite{BrSm}). 
After submitting this paper on arXiv, I received a preliminary 
version of the paper \cite{BQS} from Bridgeland, which also proves 
the result for the case $n=2$ and $2 \le N \le \infty$ 
by using a different method. 
Their approach is based on \cite{Suther}, 
while my approach is based on \cite{BrSm}.

The result is compatible with the observation in \cite{Br5} 
that if there is a Frobenius structure on the space of 
stability conditions on a $\CY_N$ category, then 
the central charge corresponds to the twisted period 
of the Frobenius structure 
with the twist parameter $\nu=(N-2) \slash 2$. 
The twisted periods 
are characterized as the horizontal sections 
of the second structure connection of the Frobenius structure 
with some parameter $\nu$ (see \cite[Section 3]{Dub1}). 
By K. Saito's theory of 
primitive forms \cite{Sai}, there is a Frobenius structure 
on $M_n$ and the map $\mathcal{W}_N$ naturally arises 
as the period map associated to the primitive form. 
As in \cite[Section 5]{Dub1}, the map $\mathcal{W}_N$ 
coincides with the twisted period map with the parameter 
$\nu=(N-2) \slash 2$. Therefore, the result supports the 
Bridgeland's conjecture in \cite{Br5}.

\subsection{Basic ideas}

Many arguments of this paper are parallel to 
those of \cite{BrSm} with appropriate generalizations. 
We replace $\CY_3$ categories with $\CY_N$ categories, 
quivers with 
$(N-2)$-colored quivers, 
simple zeros with zeros of order $N-2$ 
and triangulations with $N$-angulations. 

However as in Section \ref{sec_6}, 
there are some points where parallel arguments 
can not be applied. 
In the case of $N$-angulations, 
it is not known what objects correspond
to quivers with potential and their mutation 
in the case of triangulations.   
Hence we do not have natural dg algebras 
associated to $N$-angulations of surfaces 
and the categorical equivalence of 
Keller-Yang \cite{KY} is not established yet. 

Fortunately, King-Qiu's result \cite{KQ} on 
the hearts of $\D_n^N$ together with 
the geometric realization of the 
$(N-2)$-cluster category of type $A_n$ 
by Baur-Marsh \cite{BaMa} allows us to apply the framework of \cite{BrSm} to our cases.

The key idea of constructing stability conditions on $\D_n^N$ 
from quadratic differentials 
is the fact that saddle-free quadratic differentials 
determine $N$-angulations of surfaces as follows.  
The differential $\phi$ determines a foliation on $\P^1$ by 
\begin{equation*}
\Im\,\int p_n(z)^{\frac{N-2}{2}}dz= \text{constant}.
\end{equation*}
A leaf which connects zeros of $\phi$ is called a 
saddle trajectory. 
Near a zero of order $N-2$, the foliation 
behaves as in Figure \ref{Zero.fig}. 
As a result, 
the foliation of a saddle-free differential 
$\phi$ determines the $N$-angulation $\Delta_{\phi}$ of the  
$((N-2)(n+1)+2)$-gon (Lemma \ref{N-angulation}) 
similar to the WKB triangulations in \cite{BrSm,GMN}. 
This is the reason why we consider zeros of order $N-2$. 
Using the results of Baur-Marsh and King-Qiu, we have 
a natural heart $\H(\Delta_{\phi}) \subset \D_n^N$ 
with simple objects $S_1,\dots,S_n$ 
(Corollary \ref{correspondence}). Corresponding to 
these simple objects, there are cycles 
$\gamma_1,\dots,\gamma_n \in H_{\pm}(\phi)$, 
called standard saddle classes 
(Corollary \ref{cor_PC}). The central charge of $S_i$ 
is defined by the period of the corresponding saddle class 
$\gamma_i$:
\begin{equation*}
Z(S_i):=Z_{\phi}(\gamma_i)=\int_{\gamma_i} p_n(z)^{\frac{N-2}{2}}dz. 
\end{equation*}
Thus, we can construct 
the stability condition $(Z,\H(\Delta_{\phi}))$ 
on $\D_n^N$ (Lemma \ref{construct_stability}).

\subsection{Homotopy type problem}
Ishii-Ueda-Uehara \cite{IUU} proved that the space 
$\Stab(\D_n^2)$ is connected, 
so $\Stab^{\circ}(\D_n^2)=\Stab(\D_n^2)$. 
Together with the result of Thomas in \cite{T}, 
we have $\Stab(\D_n^2) \cong \widetilde{M_n}$. In particular, 
$\Stab(\D_n^2)$ is contractible. 
For $N \ge 3$, 
Qiu proved that $\Stab^{\circ}(\D_n^N)$ is 
simply connected (Corollary 5.5 in \cite{Q}). 
Our Theorem 1.1 
implies the following stronger result 
(Conjecture 5.8 in \cite{Q}). 
\begin{cor}[Corollary \ref{contractible}]
The distinguished connected component 
 $\Stab^{\circ}(\D_n^N)$ is contractible.
\end{cor}  
The remaining problem is the connectedness of 
the space $\Stab(\D_n^N)$. 
As in the case $N=2$, we naively expect the following. 
\begin{conjecture}
The space $\Stab(\D_n^N)$ is connected.
\end{conjecture}

\subsection*{Acknowledgements}
First of all, I am deeply grateful to my advisor Akishi Kato for 
many valuable discussions, comments, correction of errors 
and constant 
encouragement. I would like to thank Kyoji Saito and 
Atsushi Takahashi 
for valuable comments on primitive forms. I also would like to 
thank Fumihiko Sanda for helpful discussions about 
Fukaya categories. 

This work is supported by World Premier International
Research Center Initiative (WPI initiative), MEXT, Japan.

This paper is dedicated to the memory of Kentaro Nagao, 
who explained many basics of this research area 
and gave me helpful comments at 
the starting point of this work.

\subsection*{Notations}
We put $\d:=(N-2)(n+1)+2$. 
We work over an algebraically closed field $k$. 
We use homological gradings 
$d \colon M_i \to M_{i-1}$ for dg algebras and dg modules. 
Assume that all triangulated categories in this paper are  $k$-linear and 
their Grothendieck groups are free of finite rank 
($K(\D) \cong \Z^n$ for some $n$).

\section{Quadratic differentials associated to polynomials}
\label{sec_2}
In the following, we assume that $N \ge 3$ and $n \ge 1$. 
\subsection{Braid groups and their representations}
In this section, we introduce the braid groups and 
their representations on some lattices associated 
to $A_n$-quivers. 
\label{sec_braid}
\begin{defi}
The Artin braid group $B_{n+1}$ is the group generated by 
generators $\sigma_1,\sigma_2,\dots,\sigma_n$ 
with the following relations:
\begin{align*}
\sigma_i \sigma_j \sigma_i &= \sigma_j \sigma_i \sigma_j 
\quad\text{if } \abv{i-j}=1 \\
\sigma_i \sigma_j &= \sigma_i \sigma_j \quad\quad\text{if } \abv{i-j}\ge 2.
\end{align*}
\end{defi}
It is known that the element 
\begin{equation*}
C:=(\sigma_1 \sigma_2 \cdots \sigma_n )^{n+1}
\end{equation*}
generates the center of $B_{n+1}$.  

We consider two representations $\rho_+$ and $\rho_-$ 
of $B_{n+1}$ on some
lattices associated to the $A_n$-quivers. 
Such representations are equivalent to the reduced 
Burau representations with the 
parameters $\pm1$ (see \cite[Section 3]{KT}).

An $A_n$-quiver $\QA$ is a quiver obtained 
by giving orientations on edges of the Dynkin diagram 
of type $A_n$. 
We set vertices of $\QA$ by $\{1,\dots,n\}$ and 
assume that vertices $i$ and $i+1$ are connected 
by some arrow. 
Let $q_{ij}$ be the number of arrows from $i$ to $j$.

Let $L$ be a free abelian group of rank $n$ 
with generators $\alp_1,\dots,\alp_n$:
\begin{equation*}
L:=\bigoplus_{i=1}^n\,\Z \alp_i.
\end{equation*}
We define a symmetric bilinear form $\left<\,,\,\right>_+$ and 
a skew-symmetric bilinear form $\left<\,,\,\right>_-$ on $L$ by
\begin{align*}
\left<\alp_i,\alp_j \right>_+&:=\delta_{ij} +  \delta_{ji} -(q_{ij}+ q_{ji})  \\
\left<\alp_i,\alp_j \right>_-&:=\delta_{ij} - \delta_{ji} -(q_{ij}- q_{ji}) 
\end{align*} 
where $\delta_{ij}$ is the Kronecker delta. 

Reflections (which depend on the sign $\pm$)
$r_1^{\pm},\dots,r_n^{\pm} \colon L \to L$ 
via the bilinear forms $\left<\,,\,\right>_{\pm}$ is defined by 
\begin{equation*}
r_i^{\pm} (\alp_j):=\alp_j-\left<\alp_i,\alp_j \right>_{\pm}\alp_i. 
\end{equation*}
Let $W_{\pm}:=\left<r_1^{\pm},\dots,r_n^{\pm} \right>$ 
be the groups generated by these reflections. 

Two representations 
$\rho_{\pm} \colon B_{n+1} \to \GL(L,\Z)$ are defined by
\begin{equation*}
\rho_{\pm}(\sigma_i):= r_i^{\pm}. 
\end{equation*}
We put the kernel of $\rho_{\pm}$ by $P_{\pm}$. Then we have a short exact sequence
\begin{equation*}
1 \lto P_{\pm} \lto B_{n+1} \lto W_{\pm} \lto 1.
\end{equation*}

Note that the group $W_+$ is a symmetric group of degree $n+1$ and the group $P_+$ is 
a pure braid group. 

\subsection{Spaces of polynomials}
\label{sec_polynomial}
Let $p_n(z)$ be a polynomial of the following form
\begin{equation*}
p_n(z) = z^{n+1}+u_1 z^{n-1} + u_2 z^{n-2} + \cdots + u_n, \quad u_1,\dots,u_n \in \C.
\end{equation*}
We denote by $\Delta(p_n)$ the discriminant of $p_n$. (If 
$p_n(z) = \prod_{i=1}^{n+1} (z - a_i)$, then 
$\Delta(p_n) = \prod_{i < j}(a_i - a_j)^2$.)

Denote by $M_n$ the space of such polynomials with simple zeros:
\begin{equation*}
M_n := \{\,p_n(z)\,\vert\, \Delta(p_n) \neq 0\,\}.
\end{equation*}
Define the free $\C^*$-action on $M_n$ by
\begin{equation*}
(k \cdot p_n)(z):=z^{n+1}+ k u_1 z^{n-1}+ k^2 u_2 z^{n-2}+ \cdots +k^n u_n \quad(k \in \C^*).
\end{equation*}
Note that for the above action, we have the equation
\begin{equation*}
(k \cdot p_n)(k z)=k^{n+1}p_n(z).
\end{equation*}

Consider the description of $M_n$ as 
the configuration space of $\C$. 
Let $C_{n+1}(\C)$ be the configuration space of $(n+1)$ distinct points in $\C$  
and $C_{n+1}^0(\C) \subset C_{n+1}(\C)$ be the subspace consists of configurations with 
the center of mass zero:
\begin{align*}
C_{n+1}(\C) &:= \{\,(a_1,\dots,a_{n+1}) \in \C^{n+1}\,
\vert\,\, a_i \neq a_j \text{ for }i \neq j\}
\slash S_n \\
C_{n+1}^0(\C) &:= \{\,(a_1,\dots,a_{n+1}) \in \C^{n+1}\,
\vert\,\sum_{i=1}^{n+1}\,a_i=0,\, a_i \neq a_j \text{ for }i \neq j\}
\slash S_n.
\end{align*} 
A complex number $c \in \C$ acts on  $C_{n+1}(\C)$ 
as a parallel transportation 
$(a_1,\dots,a_{n+1}) \mapsto (a_1+c,\dots,a_{n+1}+c ) $.
Hence we have 
$C_{n+1}^0(\C)= C_{n+1}(\C) \slash \C$. 

As a result, we have an isomorphism 
$C_{n+1}^0(\C) \iso M_n$ given by
\begin{equation*}
(a_1,\dots,a_{n+1}) \mapsto \prod_{i=1}^{n+1}\,(z - a_i).
\end{equation*}

\begin{figure}[!ht]
\centering
\includegraphics[scale=0.9]{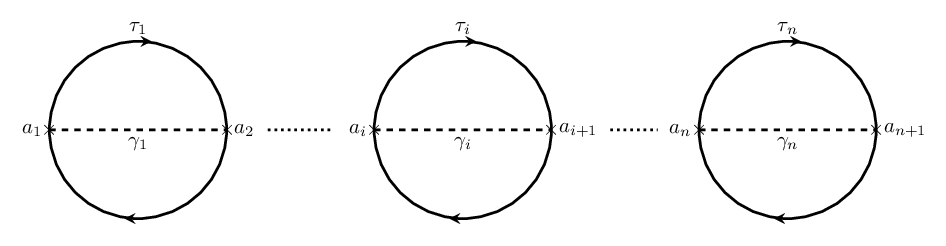}
\caption{Half-twists along the $A_n$-chain.}
\label{braid.fig}
\end{figure}

Next we consider the description of a braid group $B_{n+1}$ 
as the fundamental group of $M_n$. 
Take a configuration $\{a_1,\dots,a_{n+1}\} \in C_{n+1}(\C)$. 
We fix some order $(a_1,\dots,a_{n+1})$ for these $n+1$ points. 
Consider the collection of non-crossing paths $\gamma_1,\dots,\gamma_n$ in $\C$ such that 
$\gamma_i$ connects $a_i$ and $a_{i+1}$. We call it the $A_n$-chain of paths. 
Let $\tau_i$ be a closed path in $C_{n+1}(\C)$ 
which corresponds to the 
half-twist along $\gamma_i$ in the clockwise direction 
(see Figure \ref{braid.fig}), 
and consider its image in $M_n$ via 
the projection $C_{n+1}(\C) \to M_n$ (we also write it by $\tau_i$). 
 
It is well-known that there is an isomorphism of groups 
\begin{equation*}
\pi_1(M_n) \iso B_{n+1}, \quad [\tau_i] \mapsto \sigma_i.
\end{equation*}
For more details, we refer to \cite[Section 1.6]{KT}.

\subsection{Quadratic differentials associated to polynomials}
\label{Quadratic}
Let $\P^1$ be the Riemann sphere, and we fix the coordinate   
$z \in \C \cup \{\infty\} \cong \P^1$. 
For a polynomial $p_n \in M_n$, define the meromorphic quadratic differential $\phi$ on $\P^1$ by
\begin{equation*}
\phi(z) := p_n(z)^{N-2} dz^{\otimes 2}.
\end{equation*}
$\phi$ has zeros of order $(N-2)$ at distinct $n+1$ points on $\C \subset \P^1$, 
and a unique pole of order $(\d+2)$ at $\infty \in \P^1$.

We denote by $\Q$ the space of such quadratic differentials associated to $p_n \in M_n$:
\begin{equation*}
\Q :=\{\,p_n(z)^{N-2} dz^{\otimes 2}\,\vert\,p_n \in M_n\,\}.
\end{equation*}
Independent of $N$, the space $\Q$ is isomorphic to $M_n$ as a complex manifold, 
but we give this space various structures which depend on $N$.

We introduce the action of $\C^*$ on $\Q$. For $\phi(z)=p_n(z)^{N-2}dz^{\otimes 2}$, the action of $\C^*$ is defined by
\begin{equation*}
(k \cdot \phi)(z):=\{(k \cdot p_n)(z)\}^{N-2} dz^{\otimes 2} \quad (k \in \C^*)
\end{equation*}
where $k \cdot p_n$ is the action 
defined in Section \ref{sec_polynomial}.  
By using the equation $(k \cdot p_n)(kz)=k^{n+1}p_n(z)$, we have the following formula. 

\begin{lem}
\label{k-action}
For $\phi \in \Q$, the equation
\begin{equation*}
(k \cdot \phi)(kz)=k^{\d}\phi(z) \quad (k \in \C^*)
\end{equation*}
holds.
\end{lem}

For $\phi \in \Q$, denote by $\Zero(\phi) \subset \P^1$ the set of zeros of $\phi$ and 
by $\Crit(\phi) \subset \P^1$ the set of critical points of $\phi$. 
Note that $\Crit(\phi) = \Zero(\phi) \cup \{\infty\}$. 

The action $z \mapsto kz$ of $k \in \C^*$ on $\P^1$
 maps the critical points of $\phi$ to those of $k \cdot \phi$:
\begin{equation*}
k \,(\Crit(\phi))=\Crit(k \cdot \phi).
\end{equation*}

\subsection{Homology groups and periods}
\label{sec_homology}
\begin{defi}
Depending on the parity of $N$, we introduce homology groups $H_+(\phi)$ and $H_-(\phi)$ 
in the following. 
\begin{itemize}
\item If $N$ is even, $H_+(\phi)$ is defined by the relative homology group 
\begin{equation*}
H_+(\phi):=H_1(\C,\Zero(\phi);\Z).
\end{equation*}
\item If $N$ is odd, $H_-(\phi)$ is given as follows. 
Let $S$ be a hyperelliptic curve defined by the equation $y^2 = p_n(z)$ and 
$\pi \colon S \to \P^1$ be an associated branched covering map of degree $2$.  
Then the homology group $H_-(\phi)$ is defined by 
\begin{equation*}
H_-(\phi):= H_1(S \,\bs\, \pi^{-1}(\infty);\Z).
\end{equation*} 
\end{itemize}
\end{defi} 
Note that if $N=3$, the homology group
$H_-(\phi)$ corresponds to the hat-homology group in \cite{BrSm} which is written as $\hat{H}(\phi)$. 

\begin{lem}
Both $H_+(\phi)$ and $H_-(\phi)$ are free abelian groups of rank $n$. 
\end{lem}
{\bf Proof.} If $N$ is even, $|\Zero(\phi)|=n+1$ implies 
$\rank H_{+}(\phi)=n$.

If $N$ is odd, there are two cases, 
$n$ is even or odd. If $n$ is even, the genus of the hyperelliptic curve $S$ is $n \slash 2$ and 
the fiber $\pi^{-1}(\infty)$ consists of one point. 
Hence we have $\rank H_-(\phi)=2 ( n \slash 2)=n$. 
When $n$ is odd, 
the genus of $S$ is $(n-1) \slash 2$ and the fiber $\pi^{-1}(\infty)$ consists of two points. 
Hence we have $\rank H_-(\phi)= 2 \{ (n-1) \slash 2 \} + 1=n$. 
\endproof

Let $\psi :=\sqrt{\phi} =  p_n(z)^{\frac{N-2}{2}}\,dz$ be the square root of $\phi$. 
(If $N$ is odd, $\psi$ is determined up to the sign.)  
If $N$ is even, then $\psi$ is a holomorphic $1$-form on $\P^1 \,\bs\,\{\infty\}$, and 
if $N$ is odd, then $\psi$ is a holomorphic $1$-form on $S \,\bs\, \pi^{-1}(\infty)$. 

Therefore, for any cycle $\alp \in H_{\pm}(\phi)$ 
the integration of $1$-form $\psi$ via the cycle $\alp$ 
\begin{equation*}
Z_{\phi}(\alp) := \int_{\alp}\,\psi \in \C
\end{equation*} 
is well-defined.

Thus we have a group homomorphism
$Z_{\phi} \colon H_{\pm}(\phi) \to \C$, called the period of $\phi$.

\subsection{Framings}
Following \cite[Section 2.6]{BrSm}, we introduce 
framings on $\Q$. 
Let $\Gamma$ be a free abelian group of rank $n$. 
A $\Gamma$-framing of $\phi \in \Q$ is an isomorphism of abelian groups 
\begin{equation*}
\theta \colon \Gamma \iso H_{\pm}(\phi).
\end{equation*}
Denote by $\Q^{\Gamma}$ the space of pairs $(\phi,\theta)$ 
consisting of a quadratic 
differential $\phi \in \Q$ and a $\Gamma$-framing $\theta$. 
The projection 
\begin{equation*}
\Q^{\Gamma} \lto \Q\,,\quad \,(\phi,\theta) \mapsto \phi
\end{equation*}
defines a principal $\GL(n,\Z)$-bundle on $\Q$.  

For $(\phi,\theta) \in \Q^{\Gamma} $, 
the composition of a framing $\theta $ 
and a period $Z_{\phi} $ gives a group 
homomorphism $Z_{\phi} \circ \theta \colon \Gamma \to \C$. 
Thus, we have the map
\begin{equation*}
\mathcal{W}_N \colon \Q^{\Gamma} \lto \Hom_{\Z}(\Gamma,\C)
\end{equation*}
defined by $\mathcal{W}_N((\phi,\theta)) := Z_{\phi} \circ \theta \in \Hom_{\Z}(\Gamma,\C)$. 
We call $\mathcal{W}_N$ a period map.

\begin{prop}
\label{local_iso_period}
The period map
\begin{equation*}
\mathcal{W}_N \colon \Q^{\Gamma} \lto \Hom_{\Z}(\Gamma,\C)
\end{equation*}
is a local isomorphism of complex manifolds. 
\end{prop}
{\bf Proof. } Theorem 4.12 in \cite{BrSm} can be easily 
extended to our cases. 
In our cases, the result also follows from K. Saito's 
theory on primitive forms \cite[Section5]{Sai}
since $\mathcal{W}_N$ is interpreted as the period map 
associated to the primitive form of the $A_n$-singularity. 
\endproof

\subsection{Intersection forms}
\label{sec_intersection_form}
In this section, we introduce intersection forms on $H_{\pm}(\phi)$. 

First we consider the case $N$ is even. The exact sequence of 
relative homology groups gives the isomorphism of homology groups
\begin{equation*}
H_+(\phi) = H_1(\C,\Zero(\phi);\Z) \cong \widetilde{H}_0(\Zero(\phi);\Z)
\end{equation*}
where $\widetilde{H}_0(\Zero(\phi);\Z)$ is the reduced homology group of  
$H_0(\Zero(\phi);\Z)$, and this isomorphism is induced by the boundary map. 
There is a trivial intersection form on $H_0(\Zero(\phi);\Z)$, and it induces the 
intersection form on $\widetilde{H}_0(\Zero(\phi);\Z)$. 
Through the above isomorphism of homology groups, we have the intersection form 
\begin{equation*}
I_+ \colon H_+(\phi) \times H_+(\phi) \lto \Z.
\end{equation*}
Note that $I_+$ is symmetric and non-degenerate. 

Next consider the case $N$ is odd. Let $\pi \colon S \to \P^1$ be the hyperelliptic curve 
introduced in Section \ref{sec_homology}, and $\iota \colon S \,\bs\, \pi^{-1}(\infty) \to S$ 
be the inclusion. 
This inclusion induces the surjective group homomorphism
\begin{equation*}
\iota_* \colon H_-(\phi) = H_1( S \,\bs\, \pi^{-1}(\infty);\Z) \lto H_1(S;\Z). 
\end{equation*}
Since there is a skew-symmetric non-degenerate intersection form on $H_1(S;\Z)$, 
through the inclusion $\iota_*$, we have the intersection form
\begin{equation*}
I_- \colon H_-(\phi) \times H_-(\phi) \lto \Z.
\end{equation*}
The map $\iota_*$ is an isomorphism if $n$ is even 
and has the kernel of rank one if $n$ is odd. 
Therefore, $I_-$ is non-degenerate if $n$ is even 
or has the kernel of rank one if $n$ is odd. 

For $\phi \in \Q$, fix the order of $\Zero(\phi)$ and consider 
the $A_n$-chain of path $\gamma_1,\dots,\gamma_n$. 
If $N$ is even, the path $\gamma_i$ determines the homology 
class $\hat{\gamma}_i \in H_{+}(\phi)$ since $\gamma_i$ 
connects zeros of $\phi$. 
If $N$ is odd, since the lift $\pi^{-1}(\gamma_i)$ on the hyperelliptic curve $S$ is a closed path in $S \bs \pi^{-1}(\infty)$,  
we have the homology class 
$\hat{\gamma}_i:=[\pi^{-1}(\gamma_i)] \in H_-(\phi)$.

The pair of the homology group $H_{\pm}(\phi)$ and the intersection form $I_{\pm}$ gives 
a geometric realization of 
the lattice $(L, \left<\,,\,\right>_{\pm})$ which is introduced in Section \ref{sec_braid}.
\begin{lem}
\label{lattice_iso}
There is an isomorphism $H_{\pm}(\phi) \cong L$ given by 
$\hat{\gamma}_i \mapsto \alp_i$. 
Further, by choosing suitable orientations of $\hat{\gamma}_1,\dots,\hat{\gamma}_n$, 
this isomorphism takes the intersection form $I_{\pm}$ to the bilinear form $\left<\,,\,\right>_{\pm}$. 
\end{lem}
{\bf Proof.} The result follows from the direct calculation. 
\endproof

\subsection{Local system}
\label{sec_local_system}
For each $\phi \in \Q$, by corresponding to the homology group $H_{\pm}(\phi)$, we have 
a local system 
\begin{equation*}
\bigcup_{\phi \in \Q}H_{\pm}(\phi) \lto \Q.
\end{equation*}
For a continuous path $c \colon [0,1] \to \Q$ with  
$\phi_0=c(0),\,\phi_1=c(1)$, a parallel transportation along $c$ gives the identification  
\begin{equation*}
\GM_c \colon  H_{\pm}(\phi_0) \xto{\sim} H_{\pm}(\phi_1),
\end{equation*}
and this identification determines a flat connection on the 
local system. We call it the Gauss-Manin connection. 

In the following, we fix a point $(\phi_0,\theta_0) \in \Q^{\Gamma}$ and denote by $\Q^{\Gamma}_*$ 
the connected component containing $(\phi_0,\theta_0)$. 
Any point $(\phi,\theta) \in \Q^{\Gamma}$ is contained in $\Q^{\Gamma}_*$ 
if and only if there is a continuous path $c \colon [0,1] \to \Q$ connecting 
$\phi_0$ and $\phi$ such that the diagram  
\begin{equation*}
\xymatrix{ 
& \Gamma \ar[dl]_{\theta_0} \ar[dr]^{\theta}  & \\ 
H_{\pm}(\phi_0) \ar[rr]^{\GM_c} && H_{\pm}(\phi)
}
\end{equation*}
commutes. 

Note that different connected components of $\Q^{\Gamma}$ are all isomorphic since 
these are related by the action of $\GL(n,\Z)$. 

We introduce the $\C$-action on $\Q^{\Gamma}_*$, which is compatible with the $\C$-action 
on the space of stability conditions given in Section \ref{sec_group}. 
\begin{defi}
\label{C_action}
For a framed differential $(\phi,\theta) \in \Q^{\Gamma}_*$, the new framed differential 
$(\phi^{\prime},\theta^{\prime}):=t \cdot (\phi,\theta)$ for 
$t \in \C$ is defined by 
\begin{equation*}
\phi^{\prime}:=e^{-(2 \pi i\,\slash \d)\, t}\cdot \phi,\quad \theta^{\prime}:=
\GM_c \circ \theta
\end{equation*}
where $\GM_c$ is the Gauss-Manin connection along the path 
$c(s):=e^{-(2 \pi i\,\slash \d)\,t\,s} \cdot \phi$ for $s \in [0,1]$ which connects $\phi$ and $\phi^{\prime}$. 
\end{defi}

\subsection{Monodromy representations}
By using the intersection form $I_{\pm}$ on $H_{\pm}(\phi)$ defined in Section \ref{sec_intersection_form}, 
we can compute how cycles 
change by the Gauss-Manin connection 
along a closed path in $\Q$. The result is 
known as the Picard-Lefschetz formula. 

Let $\phi_0=\left[\prod_{i=1}^{n+1}(z-a_i)\right]^{N-2} dz^{\otimes 2}$ be the fixed differential in the previous section. 
Take the $A_n$-chain of paths  $\gamma_1,\dots,\gamma_n$ 
for ordered zeros $(a_1,a_2,\dots,a_{n+1})$
as 
in Section \ref{sec_polynomial}.
Under the identification $M_n \cong \Q$, 
we consider the half-twists $\tau_1,\dots,
\tau_n$ as in Section \ref{sec_polynomial} 
and take the orientations of $\hat{\gamma}_1,\dots,\hat{\gamma}_n \in H_{\pm}(\phi_0)$ as in 
Lemma \ref{lattice_iso}.

\begin{prop}
\label{monodromy}
The monodromy representation of $[\tau_i] 
\in \pi_1(\Q,\phi_0)$ on $H_{\pm}(\phi_0)$ is given by
\begin{equation*}
\GM_{\tau_i}(\beta) = \beta - I_{\pm}(\hat{\gamma}_i,\beta)\,\hat{\gamma}_i \quad(\,\beta \in H_{\pm}(\phi_0)\,).
\end{equation*}
\end{prop}

Lemma \ref{lattice_iso} and Proposition \ref{monodromy} imply that 
the monodromy representation of the braid group $\pi_1(M_n,*) \cong B_{n+1}$ on $H_{\pm}(\phi_0)$ gives 
a geometric realization of two representations $\rho_+$ 
and $\rho_-$ introduced in 
Section \ref{sec_braid}. Further recall from Section \ref{sec_braid} that the subgroup $P_{\pm} \subset B_{n+1}$ is 
the kernel of $\rho_{\pm}$ and the group $W_{\pm}$ is the image 
of $\rho_{\pm}$ on $L \cong H_{\pm}(\phi_0)$.

Let $\widetilde{M_n}$ be a universal cover of $M_n$. 
Note that $B_{n+1}$ acts on $\widetilde{M_n}$ as the group of 
deck transformations. 

\begin{cor}
\label{universal}
The connected component $\Q_*^{\Gamma}$ is isomorphic to the space 
$\widetilde{M_n} \slash P_{\pm}$ and 
the projection map $\Q_*^{\Gamma} \to \Q$ gives a principal $W_{\pm}$-bundle. 
\end{cor}
{\bf Proof.} It follows from the above remarks. 
\endproof

\section{Trajectories on $\P^1$}
In Section 3, in order to fix notations, we collect basic facts and definitions for trajectories on $\P^1$ from 
\cite[Section 3]{BrSm}. 
Since our Riemann surface is only $\P^1$, many results are simplified.   

\subsection{Foliations}
Let $\phi \in \Q$. For any neighborhood of a point in $\P^1 \, \bs\, \Crit (\phi)$, 
the differential $\phi(z)=\varphi(z)dz^{\otimes 2}$ defines the distinguished local coordinate $w$ by
\begin{equation*}
w := \int \sqrt{\varphi(z)}\,dz,
\end{equation*}
up to the transformation $w \mapsto \pm w +c\,(c \in \C)$. 
The coordinate $w$ also characterized by  
\begin{equation*}
\phi(w)= dw \otimes dw.
\end{equation*} 

The coordinate $w$ determines a foliation 
of the phase $\theta$ on $\P^1 \, \bs\, \Crit (\phi)$ by 
\begin{equation*}
\Im (w \slash e^{i \pi \theta})= \text{constant}.
\end{equation*}

A straight arc $\gamma \colon I \to \P^1 \,\bs\,\Crit(\phi) $ of the phase $\theta$ 
is an integral curve along the foliation $\Im (w \slash e^{i \pi \theta})= \text{constant}$, 
defined on an open interval $I \subset \R$.  
Note that a maximal  straight arc of the phase $\theta$ is 
a leaf of the foliation of the phase $\theta$.   

A foliation of the phase $\theta = 0$ is called 
a horizontal foliation, 
and a maximal straight arc of the phase $\theta=0$ is called a trajectory.

\begin{lem}
\label{rotation}
Let $\phi \in \Q $. For $k=r e^{i \pi \theta} \in \C^*$, set 
$\phi^{\prime}:=k \cdot \phi$. 
Then the distinguished coordinates $w$ and $w^{\prime}$ which 
come from $\phi$ and $\phi^{\prime}$ 
have the relation
\begin{equation*}
w^{\prime}=k^{\frac{\d}{2}} w =r^{\frac{\d }{2}}e^{i \pi  \frac{\d \theta }{2}} w.
\end{equation*}
In addition, the transformation $z \mapsto k z$ on $\P^1$ 
maps the horizontal straight arcs for $\phi$ to 
the straight arcs of the phase $(\d \theta )\slash 2$ for $\phi^{\prime}$. 
\end{lem}
{\bf Proof.} 
Let
$\gamma \colon I \to \P^1 \,\bs\,\Crit(\phi) $ be a horizontal straight arc of $\phi$. 
Then the following computation gives the result:
\begin{align*}
\int^{k \gamma(t)} \sqrt{\varphi^{\prime}(z^{\prime})}\,dz^{\prime} 
&= \int^{k \gamma(t)} \sqrt{(k \cdot \varphi)
(z^{\prime})}\,dz^{\prime} \\
&=\int^{\gamma(t)} \sqrt{(k \cdot \varphi)(k z)}\,
d(k z)  \quad (\text{put } z^{\prime}=kz.) \\ 
&= k^{\frac{\d}{2}} \int^{\gamma(t)}  \varphi(z)\,d z 
\quad (\text{by Lemma \ref{k-action}}).
\end{align*}
\endproof

A maximal straight arc $\gamma$ of the phase $\theta$ defined 
on a finite open interval $(a,b)$ can be extended 
to a continuous path $\gamma \colon [a,b] \to \P^1$ with $\gamma(a),\gamma(b) \in \Zero(\phi)$ since $\P^1$ is compact. 
According to $\theta=0$ or $\theta \neq 0$, 
we call $\gamma$ a saddle trajectory or a saddle connection 
respectively.

\subsection{Foliations near critical points}
\label{critical}
\begin{figure}[!ht]
\centering
\includegraphics[scale=0.9]{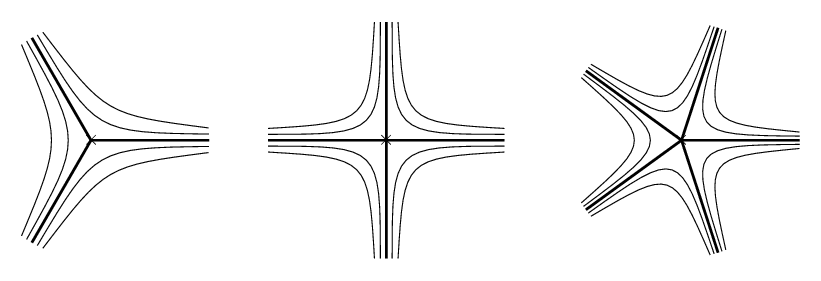}
\caption{Trajectories near zeros of order $1,2,3.$}
\label{Zero.fig}
\end{figure}

Here we see the local behavior of the horizontal foliation near a critical point. In details, we refer to \cite[Section 3.3]{BrSm}.
First we consider the behavior near a zero. 
Let $\phi \in \Q$ and $p \in \Zero(\phi)$ be a zero of order $N-2$. Then 
there is a local coordinate $t$ of a neighborhood $p \in U \subset \P^1$ ($t=0$ corresponds to $p$) 
such that 
\begin{equation*}
\phi=c^2\, t^{N-2} dt^{\otimes 2}, \quad c=\frac{1}{2}\,N.
\end{equation*}
Hence on $U \,\bs\,\{ p\}$, the distinguished local coordinate $w$ takes the 
form $w=t^{N \slash 2}$. Figure \ref{Zero.fig} 
illustrates the local behavior of the horizontal foliation 
in the case $N=3,4,5$. 

Note that there are $N$ horizontal straight arcs departing from $p$. 
This fact is the reason why the differential $\phi$ determines 
 $N$-angulations. 

Next we consider the behavior near a pole. 
Recall from Section \ref{Quadratic} that 
$\phi$ has a unique pole of order $\d+2$ at $\infty \in \P^1$. 
Therefore it is sufficient to consider a pole of order $>2$. 
In this case, there is a neighborhood $\infty \in U \subset \P^1$ and a collection of $\d$ distinguished 
tangent directions $v_i\,(i=1,\dots,\d)$ at $\infty$, such that any trajectory entering $U$ tends to 
$\infty$ and becomes asymptotic to one of the $v_i$. 
For more details, we refer to \cite[Section 3.3]{BrSm}.

\subsection{Classification of trajectories}
In this section, we summarize the global structure of the horizontal foliation on $\P^1$ 
following \cite[Section 3.4]{BrSm}, 
and \cite[Section 9-11 and 15]{Str}.

We first note that $\phi$ is holomorphic on 
the contractible domain $\P^1 \,\bs\,\{\infty\}$. 
Then Theorem 14.2.2 in \cite{Str} implies that 
there are no closed trajectories 
and Theorem 15.2 in \cite{Str} implies that there are 
no recurrent trajectories. 
Hence in our setting, every trajectory of $\phi$ is classified to 
exactly one of the following classes. 
\begin{itemize}
\item[(1)] A saddle trajectory approaches distinct zeros of $\phi$ at both ends. 
\item[(2)] A separating trajectory approaches a zero of $\phi$ at one end and 
a pole $\infty \in \P^1$ at the other end.
\item[(3)] A generic trajectory approaches a pole $\infty \in \P^1$ at both ends. 
\end{itemize}

As explained in the previous section, there are only finite horizontal straight arcs departing from zeros of $\phi$. 
Therefore, the number of saddle trajectories and separating trajectories is finite. 
For a quadratic differential $\phi$, denote by $s_{\phi}$ the number of saddle trajectories, and by $t_\phi$ the 
number of separating trajectories. 
Then the fact that there are $N$ horizontal straight arcs departing from each zero of $\phi$ implies 
the equation
\begin{equation}
\label{eq3.1}
2s_\phi + t_\phi =N(n+1).
\end{equation}

By removing these two type of trajectories and  $\Crit(\phi)$ from $\P^1$, 
the remaining open surface splits into a disjoint union of connected components. 
In our setting, each connected component coincides with one of the following two type of surfaces. 
\begin{itemize}
\item[(1)] A half-plane is the upper half-plane 
\begin{equation*}
\{\,w \in \C\,\vert\,\Im w >0\,\}
\end{equation*}
equipped with the differential $dw^{\otimes 2}$. It is swept out by generic trajectories 
with the both end points at $\infty \in \P^1 $. The boundary 
 $\{\Im w =0\}$ consists of saddle trajectories and two separating trajectories. 
\item[(2)] A horizontal strip is the strip domain
\begin{equation*}
\{\,w \in \C\,\vert\,a <\Im w <b \,\}
\end{equation*}
equipped with the differential $dw^{\otimes 2}$. It is also swept out 
by generic trajectories 
with the both end points at $\infty \in \P^1 $. Each of the two boundaries 
$\{\Im w=a\}$ and $\{\Im w=b \}$ consists of 
saddle trajectories and two separating trajectories. 
\end{itemize}
Note that every saddle trajectory or separating trajectory appears just two times 
as a boundary element of these open surfaces.    

\begin{figure}[!ht]
\centering
\includegraphics[scale=0.9]{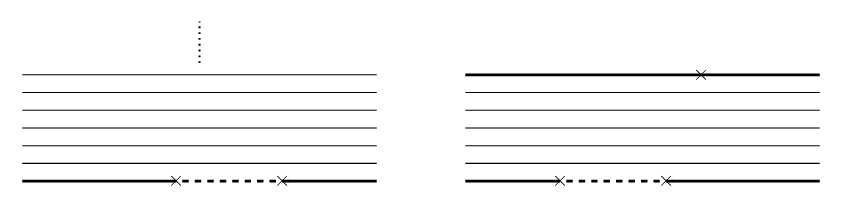}
\caption{Examples of a half plane and a horizontal strip.}
\label{domain.fig}
\end{figure}

Let $k$ be the number of half planes and $l$ be the number of horizontal strips appearing 
in such a decomposition of $\P^1$.  
Since every half plane has $2$ separating trajectories on the boundary, and every horizontal strip 
has $4$ separating trajectories on the two boundaries, 
we have the equation
\begin{equation}
\label{eq3.2}
2k + 4l=2 t_\phi.
\end{equation}

We also note that near  a pole of order $m >2$, just $m-2$ half planes appear.  
Since $\phi$ has a unique pole of order $(\d+2)$ at $\infty \in \P^1$, we have 
\begin{equation}
\label{eq3.3}
k=\d=(N-2)(n+1)+2.
\end{equation}

\begin{lem}
\label{saddle_number}
For $\phi \in \Q$, the number of  saddle trajectories $s_\phi$ and 
the number of horizontal strips $l$ satisfy the equation
\begin{equation*}
s_\phi +l =n.
\end{equation*}
\end{lem}
{\bf Proof.} 
This immediately follows from equations (\ref{eq3.1}), (\ref{eq3.2}) and (\ref{eq3.3}).
\endproof

\begin{defi}
A differential $\phi \in \Q$ is called saddle-free if $\phi$ has no saddle trajectories (i.e. $s_\phi=0$). 
\end{defi}
For a saddle-free differential $\phi \in \Q$, Lemma \ref{saddle_number}  
implies that exactly $n$ horizontal strips appear in the decomposition of $\P^1$ by $\phi$. 
In this situation, any boundary component of a half plane or a horizontal strip consists of 
two separating trajectories. 

\subsection{Hat-homology classes and standard saddle connections}
\label{sec_standard_saddle}
In this section, we collect some constructions of homology classes from \cite[Section 3.2 and 3.4]{BrSm}. 

Take $\phi \in \Q$. Recall from Section \ref{sec_homology} that 
the square root of $\phi$ defines a holomorphic $1$-form $\psi=\sqrt{\phi}$ on 
$\P^1 \,\bs\, \{\infty\}$ or $S\,\bs\,\pi^{-1}(\infty)$ depending on the parity of $N$. 
We set 
\begin{equation*}
\hat{S}^{\circ}:=
\begin{cases}
\P^1 \,\bs\, \{\infty\} &\text{ if $N$ is even} \\
S\,\bs\,\pi^{-1}(\infty) &\text{ if $N$ is odd}.
\end{cases}
\end{equation*}

We first note that the $1$-form $\psi$ defines the distinguished 
local coordinate 
$\hat{w}$ on $\hat{S}^{\circ}\,\bs\,\Zero(\psi)$ by $\psi=d \hat{w}$, up to the transformation 
$\hat{w} \mapsto \hat{w}+c\,(c \in \C)$. As a result, 
$\psi$ determines the horizontal 
foliation on 
$\hat{S}^{\circ}\,\bs\,\Zero(\psi)$ by 
$\Im \hat{w}=\text{constant}$. 
This foliation has the canonical orientation 
to make the tangent vectors lie in $\R_{>0}$. 
We call it the oriented horizontal foliation 
of $\psi$. 

For a saddle trajectory 
$\gamma \colon [a,b] \to \P^1$, we can define a 
homology class $\hat{\gamma} \in H_{\pm}(\phi)$, called the hat-homology class of 
$\gamma$  by the following.  

If $N$ is even, we can orient 
$\gamma$ to be compatible with the 
oriented horizontal foliation of $\psi$.
Since $\gamma(a),\gamma(b) \in \Zero(\phi)$, 
the oriented curve $\gamma$ defines a homology class
\begin{equation*} 
\hat{\gamma}:=[\gamma] \in H_+(\phi)=H_1(\C,\Zero(\phi);\Z).
\end{equation*}. 
By definition, it satisfies $Z_{\phi}(\hat{\gamma}) \in \R_{>0}$. 

If $N$ is odd, since zeros of $\phi$ are branched points of the double cover $\pi \colon S \to \P^1$,  
the inverse image $\pi^{-1}(\gamma)$ becomes a closed curve in 
$\hat{S}^{\circ}$. 
We can orient $\pi^{-1}(\gamma)$ to be compatible with 
the oriented horizontal orientation of $\psi$. Then 
the oriented closed curve $\pi^{-1}(\gamma)$ 
defines a homology class
\begin{equation*}
\hat{\gamma}:=[\pi^{-1}(\gamma)] 
\in H_-(\phi)=H_1(S \,\bs \,\pi^{-1}(\infty);\Z). 
\end{equation*}
In the case that $N$ is even, it satisfies $Z_{\phi}(\hat{\gamma}) \in \R_{>0}$. 

A similar construction works for a saddle connection $\gamma \colon [a,b] \to \P^1$ with 
the phase $\theta \in (0,\pi)$. In this case, 
the orientation is determined to satisfy that 
$Z_{\phi}(\hat{\gamma}) \in \{\,r e^{i \pi \theta}\,\vert\,r >0\,\}\subset \mathbb{H}$.  

Let $h=\{ a <\Im w < b \}$ be a horizontal strip appearing in 
the domain decomposition of $\P^1$ by 
a saddle-free differential $\phi$. 
Then each of two boundary components $\{\Im w=a\}$ and $\{\Im w=b\}$ contains 
exactly one zero of $\phi$ and consists of two separating trajectories intersecting 
at the zero. 
Hence, there is a unique saddle connection $l_h$ with the phase  $\theta \in (0,\pi) $ which connects two zeros appearing in two boundary components. 
We denote by $\gamma_h:=\hat{l_h} \in H_{\pm}(\phi)$  
the hat-homology class associated to the saddle connection $l_h$. 

\begin{figure}[!ht]
\centering
\includegraphics[scale=0.9]{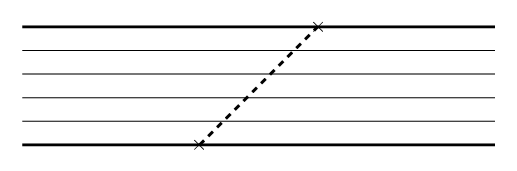}
\caption{Standard saddle connection.} \label{standard.fig}
\end{figure}

Since $\phi$ is saddle-free, Lemma \ref{saddle_number} implies that there are precisely 
$n$ horizontal strips in the decomposition of $\P^1$ by $\phi$. 
Thus, we obtain $n$ 
saddle connections $l_{h_1},\dots,l_{h_n}$ and $n$ hat-homology classes $\gamma_{h_1},\dots,\gamma_{h_n}$. 
$l_{h_i}$ is called a standard saddle connection, and 
$\gamma_{h_i}$ is called a standard saddle class.

\begin{lem}[\cite{BrSm}, Lemma 3.2]
\label{standard_saddle}
Standard saddle classes $\gamma_{h_1},\dots,\gamma_{h_n}$ form a basis of $H_{\pm}(\phi)$. 
\end{lem}
{\bf Proof. } The argument of Lemma 3.2 in \cite{BrSm} for $N=3$ also works for $N \ge 3$.
\endproof

\section{$CY_N$ categories associated to acyclic quivers}
\label{sec_4}
For objects $E,F $ in a triangulated category $\D$, 
we set
\begin{align*}
\Hom^i(E,F)&:=\Hom_{\D}(E,F[i]) \\
\dhom^i(E,F)&:=\dim_k \Hom_{\D}(E,F[i])
\end{align*} 
where $[i]$ is the $i$-th shift functor. 

\subsection{Ginzburg dg algebras}
\label{Ginzburg}
Let $Q=(Q_0,Q_1)$ be a finite acyclic quiver with 
vertices $Q_0=\{1,\dots,n\}$ and arrows $Q_1$.  
The Ginzburg dg algebra \cite{Ginz} 
$\GQ :=(k \overline{Q},d)$ is defined as follows. 
Define a graded quiver $\overline{Q}$ with vertices 
$\overline{Q}_0=\{1,\dots,n\}$ and following arrows:
\begin{itemize}
\item an original arrow $a \colon i \to j \in Q_1$ (degree $0$); 
\item an opposite arrow $a^* \colon j \to i$ for the original arrow $a \colon i \to j \in Q_1$ (degree $N-2$);
\item a loop $t_i$ for each vertex $i \in Q_0$ (degree $N-1$).
\end{itemize}
Let $k \overline{Q}$ be a graded path algebra of $\overline{Q}$, 
and define a differential $d \colon k \overline{Q} \to k \overline{Q}$ of degree $-1$ by
\begin{itemize}
\item $da =da^*= 0$ for $a \in Q_1$;
\item $d t_i = e_i \,\left(\sum_{a \in Q_1}(aa^* -a^*a)\right)\,e_i$;
\end{itemize}
where $e_i$ is the idempotent at $i \in Q_0$. Thus 
we have the dg algebra $\GQ =(k \overline{Q},d)$. 
Note that the $0$-th homology is given by 
$\Ho_0(\GQ) \cong k Q$. 

Let $\D(\GQ)$ be the derived category of right dg-modules over $\GQ$, and 
$\DQ$ be the full subcategory of $\D(\GQ)$ consists of dg modules $M$ whose 
homology is of finite total dimension:
\begin{equation*}
\DQ:=\{\, M \in \D(\GQ) \,\vert\,
\sum_{i \in \Z}\,\dim_k \Ho_i(M) <\infty \,\}.
\end{equation*}

A triangulated category $\D$ is called $N$-Calabi-Yau ($\CY_N$) if 
for any objects $E,F \in \D$, there is a natural isomorphism of $k$-vector spaces
\begin{equation*}
\nu \colon \Hom(E,F) \xto{\sim} \Hom^N(F,E)^*.
\end{equation*}
(For a vector space $V$, we denote by $V^*$ the dual vector space.)

\begin{thm}[\cite{Kel2}, Theorem 6.3]
The category $\DQ$ is a $\CY_N$ triangulated category. 
\end{thm}

Let $K(\DQ)$  be the Grothendieck group of $\DQ$.  
A bilinear form  
$\chi \colon K(\DQ) \times K(\DQ) \to \Z$, 
called the Euler form, is 
defined to be
\begin{equation*}
\chi(E,F):=\sum_{i\in \Z}(-1)^i \dhom^i(E,F), 
\quad \text{for } E,F \in \DQ.
\end{equation*}
The $\CY_N$ property of $\DQ$ implies that 
$\chi$ is symmetric if $N$ is even, and 
skew-symmetric if $N$ is odd.


For a quiver $\QA$, we set $\D_n^N:= \D^N_{\QA}$. 

\subsection{Hearts of bounded t-structures}

\begin{defi}[\cite{BBD}]
A t-structure on $\D$ is a full subcategory 
$\sli \subset \D$ satisfying the following conditions:
\begin{itemize}
\item[(a)] $\sli[1] \subset \sli$,
\item[(b)] define $\sli^{\perp} :=\{G \in \D \vert \Hom(F,G) =0 \text{ for all }F \in \sli\,  \}$, 
then for every object $E \in \D$ there is an exact triangle $F \to E \to G \to F[1]$ in $\D$ 
with $F \in \sli$ and $G \in \sli^{\perp}$.
\end{itemize}
In addition, the t-structure $\sli \subset \D$ is said to be bounded if $\sli$ satisfies the condition:
\begin{equation*}
\D = \bigcup_{i,j \in \Z} \sli^{\perp}[i] \cap \sli[j].
\end{equation*}
\end{defi}

For a t-structure $\sli \subset \D$, we define its heart $\H$ by
\begin{equation*}
\H := \sli^{\perp}[1] \cap \sli.
\end{equation*}
It was proved in \cite{BBD} that $\H$ becomes an abelian category. 

\begin{rem}
\label{iso_K-group}
For a heart $\H \subset \D$ of a bounded t-structure, 
we have a canonical isomorphism of Grothendieck groups
\begin{equation*}
K(\H) \cong K(\D).
\end{equation*}
\end{rem}

In the following, we only treat bounded t-structures and their hearts. 
The word heart always refers to the heart of some 
bounded t-structure.

\begin{prop}[\cite{Am}, Section 2]
\label{standard_heart}
There is a canonical bounded t-structure $\sli_0 \subset \DQ$ with the heart $\H_0$ such that 
the functor $\Ho_0 \colon \DQ \to \Mod k Q$ 
induces an equivalence of abelian categories:
\begin{equation*}
  \Ho_0 \colon \H_0 \iso \Mod k Q.
\end{equation*}
\end{prop}
We call $\H_0$ the standard heart of $\DQ$. 
$\H_0$ is generated by simple $\GQ$-modules 
$R_1,\dots,R_n$ which correspond to vertices 
$1,\dots,n$ of $Q$.

\subsection{Spherical twists}
\label{sec_spherical}
We define some autoequivalences of $\DQ$, called spherical twists, 
introduced by P. Seidel and R.Thomas \cite{ST}.

An object $S \in \DQ$ is called $N$-spherical if 
\begin{align*}
\Hom^i(S,S) = 
\begin{cases}
k \quad \text{if} \quad i=0,N  \\
\,0 \quad \text{otherwise} . 
\end{cases}
\end{align*}

\begin{prop}[\cite{ST}, Proposition 2.10]
For a spherical object $S \in \DQ$, there is an 
exact autoequivalence $\Phi_S \in \Aut(\DQ)$ 
defined by the exact triangle 
\begin{equation*}
\Hom^{\bullet}(S,E) \otimes S \lto E \lto \Phi_S(E)
\end{equation*}
for any object $E \in \DQ$. The inverse functor $\Phi_S^{-1} \in \Aut(\D^N_n)$ is given by
\begin{equation*}
\Phi_S^{-1}(E) \lto E \lto S \otimes \Hom^{\bullet}(E,S)^* .
\end{equation*}
\end{prop}

Let $R_1,\dots,R_n \in \H_0$ be simple $\GQ$-modules. 
Then $R_1,\dots,R_n$ 
are $N$-spherical in $\DQ$ (see Lemma 4.4 in \cite{Kel2}).  
Hence they define spherical twists $\Phi_{R_1},\dots,\Phi_{R_n} \in \Aut(\DQ)$. 
The Seidel-Thomas braid group $\Sph(\DQ)$ is defined 
to be the subgroup of $\Aut(\DQ)$ generated by these spherical twists:
\begin{equation*}
\Sph(\DQ) := \left< \Phi_{R_1},\dots,\Phi_{R_n} \right>.
\end{equation*}

In the rest of this section, we assume that $Q=\QA$. 
\begin{thm}[\cite{ST}, Theorem 1.2 and Theorem 1.3]
\label{faithful}
The correspondence of generators 
$\sigma_i \mapsto \Phi_{R_i}$ 
gives an isomorphism of groups
\begin{equation*}
B_{n+1} \cong \Sph(\D^N_n).
\end{equation*}
\end{thm}


On $K(\D^N_n)$, a spherical twist $\Phi_{R_i}$ induces the reflection $[\Phi_{R_i}] \colon K(\D^N_n) \to K(\D^N_n)$ given by
\begin{equation*}
[\Phi_{R_i}]([E])=[E]- \chi(R_i,E)[R_i].
\end{equation*}
Recall from Section \ref{sec_braid} that for a quiver $\QA$, 
we define the lattice $(L,\left<\,,\,\right>_{\pm})$. 
By the map $[R_i] \mapsto \alp_i$, we have an isomorphism 
of abelian groups 
\begin{equation*}
K(\D_n^N) \cong L.
\end{equation*}
Further, this map takes the Euler form $\chi$ to the 
bilinear form $\left<\,,\,\right>_{\pm}$. 
Therefore, under the above isomorphism, 
spherical twists $\Phi_{R_1}, \dots,\Phi_{R_n}$ act 
on $K(\D^N_n)$
as reflections $r^{\pm}_1,\dots,r^{\pm}_n$ 
in Section \ref{sec_braid} :
\begin{equation*}
[\Phi_{R_1}]=r^{\pm}_1,\dots,[\Phi_{R_n}]=r^{\pm}_n.
\end{equation*}

\subsection{Koszul duality and spherical twists}
\label{Koszul}
As in Section \ref{Ginzburg}, 
let $\GQ$ be the Ginzburg dg $\CY_N$ algebra of 
an acyclic quiver $Q$, $\D(\GQ)$ be the 
derived category of right dg $\GQ$-modules and 
$\DQ \subset \D(\GQ)$ be the full subcategory consisting 
of right dg $\GQ$-modules 
with finite dimensional total cohomology. 

In the rest of this section, we write $A:=\GQ$. 

Consider the standard heart $\H_0 \subset \DQ$ and 
its simple objects $R_1,R_2,\dots,R_n \in \H_0$. 
Set $R:=\oplus_{i=1}^n R_i$ and define the endomorphism 
dg algebra $B$ by  
\begin{equation*} 
B:=\R \Hom_A(R,R).  
\end{equation*}
(We can easily check that the dg algebra $B$ is quasi-isomorphic to the underlying  graded algebra, 
i.e. the differential of $B$ is zero.)

Denote by $\D(B)$ the derived category of right 
dg $B$-modules. 
The perfect derived category $\Per(B) $ is defined to be 
the smallest triangulated full subcategory of $\D(B)$ containing 
$B$ and closed under taking direct summands. 

The derived functor 
\begin{equation*}
F:=\R \Hom_A(R,-) \colon \D(A) \to \D(B)
\end{equation*}
has a right adjoint functor 
\begin{equation*}
G:=-\otimes_B^{\mathbb{L}}R \colon \D(B) \to \D(A).
\end{equation*}
By restricting $F$ on $\DQ \subset \D(A)$ and 
$G$ on $\Per(B) \subset \D(B)$, 
we have the following (see Section 9 in \cite{KN}). 
\begin{prop}
\label{equiv}
The restriction of $F$ on 
$\DQ \subset \D(A)$ gives  
an equivalence of triangulated categories:
 \begin{equation*} 
F \colon \DQ \iso \Per(B).
\end{equation*}
Its inverse is given by the restriction of $G$ on 
$\Per(B)\subset \D(B) $:
\begin{equation*}
G \colon \Per(B) \iso \DQ.
\end{equation*}
\end{prop}

\begin{rem}
If $Q=\QA$, then the graded algebra $B$ coincides with the 
graded algebra $A_n^N$  in \cite[Section 3]{T} 
and the category $\Per(B)$ coincides with 
the category $D^N_n$ in \cite[Definition 3.2]{T}. 
\end{rem}

We set $P_i:=F(R_i)=\R \Hom_A(R,R_i) \in \Per(B)$. 
Then we have $B=\oplus_{i=1}^n P_i$. 

For each $P_i$, define the object 
$Q_i \in \Per (B^{op} \otimes B)$ by the cone 
\begin{equation*}
Q_i:=\Cone (\,\R \Hom_B(P_i,B) \otimes P_i \to B).
\end{equation*}

Since $R_i$ is $N$-spherical in $\DQ$, the object 
$P_i$ is also $N$-spherical in $\Per(B)$. We denote by $\Psi_{P_i} \in \Aut(\Per(B))$ 
the spherical twist along $P_i$.
\begin{lem}
\label{tensor}
The spherical twist $\Psi_{P_i} $ 
is given by the derived tensor 
functor $-\otimes_B^{\mathbb{L}}Q_i$:
\begin{equation*}
\Psi_{P_i} \cong - \otimes_B^{\mathbb{L}}Q_i.
\end{equation*}
\end{lem}
{\bf Proof.} 
Since $P_i$ is perfect, there is an isomorphism 
\begin{equation*}
N \otimes_B^{\mathbb{L}} \R \Hom_B(P_i,B) 
\cong \R \Hom_B(P_i,N)
\end{equation*}
for $N \in \Per(B)$. 
Therefore we have 
\begin{align*}
N \otimes_B^{\mathbb{L}}Q_i 
&= \Cone (N \otimes_B^{\mathbb{L}} \R \Hom_B(P_i,B) 
\otimes P_i \to N\otimes_B^{\mathbb{L}} B) \\
&\cong \Cone (\,\R \Hom_B(P_i,N) \otimes P_i \to N) \\
&=\Psi_{P_i}(N). 
\end{align*}
Since the above computation is functorial for $N$, 
the result follows. 
\endproof

Define the Seidel-Thomas braid group of $\Per(B)$ by  
\begin{equation*}
\Br(\Per(B)):=\left<\Psi_{P_1},\dots,\Psi_{P_n} \right>\subset \Aut(\Per(B)).
\end{equation*}

\begin{prop}
\label{tensor_id}
If $\Psi \in \Br(\Per(B))$ preserves $P_1,\dots,P_n$, i.e. 
there is an isomorphism $\Psi(P_i) \cong P_i$ for each $P_i$, 
then $\Psi$ is  isomorphic to the identity functor : 
$\Psi \cong \id$. 
\end{prop}
{\bf Proof. }Assume that 
\begin{equation*}
\Psi=\Psi_{P_{i_1}} \circ \dots \circ \Psi_{P_{i_s}}. 
\end{equation*}
Then by Lemma \ref{tensor}, $\Psi$ is equivalent to 
the derived tensor functor $-\otimes_B^{\mathbb{L}} Q$ 
where the object $Q \in \Per(B^{op} \otimes B)$ is given by 
\begin{equation*}
Q:=Q_{i_s} \otimes_B^{\mathbb{L}} \dots \otimes_B^{\mathbb{L}} 
Q_{i_1}.
\end{equation*}
Since $\Psi(P_i) \cong P_i$, we have an isomorphism
$P_i \otimes_B^{\mathbb{L}} Q \cong P_i$ for each $P_i$. 
Recall that $B=\oplus_{i=1}^n P_i$. 
By taking a direct sum, we have 
\begin{equation*}
Q \cong B \otimes_B^{\mathbb{L}} Q=
(\oplus_{i=1}^n P_i)\otimes_B^{\mathbb{L}} Q \cong 
(\oplus_{i=1}^n P_i) =B. 
\end{equation*}
Thus, $\Psi \cong -\otimes_B^{\mathbb{L}} B$ since $Q \cong B$. For each $N \in \Per(B) $, 
\begin{equation*}
N \to N \otimes_B^{\mathbb{L}} B,\quad n \mapsto n \otimes 1 
\end{equation*}
gives a natural isomorphism $\id \cong -\otimes_B^{\mathbb{L}} B$.
\endproof

\begin{cor}
\label{tensor_id2}
If $\Phi \in \Br(\DQ)$ preserves $R_1,\dots,R_n$, i.e. 
there is an isomorphism $\Phi(R_i) \cong R_i$ for each $R_i$, 
then $\Phi$ is  isomorphic to the identity functor : 
$\Phi \cong \id$. 
\end{cor}

\subsection{Simple tilted hearts}
\label{sec_tilting}

Let $\H \subset \D$ be a heart in a triangulated category $\D$. 
We denote by $\Sim \H$ 
the set of simple objects in $\H$. 
\begin{defi}
$\H$ is called finite if $\Sim \H$ 
is a finite set and any object of $\H$ is generated by 
finitely many extensions of $\Sim \H$.  
\end{defi}

For a given heart, one can construct a new heart by 
using the technique of tilting in \cite{HRS}. 
In particular for a finite heart $\H \subset \D$ and a
simple object $S \in \H$, 
there is a new heart
$\H^{\sharp}_S$ (respectively $\H^{\flat}_S$), 
called the forward (respectively backward) simple tilt 
of $\H$ by $S$. 
For the precise definition of $\H^{\sharp}_S$ and $\H^{\flat}_S$,  
we refer to \cite[Section 3]{KQ}. (Our notation is 
based on \cite{KQ}.)  
In the present paper, we do not use the explicit construction 
of simple tilted hearts. 

We say an object $S \in \H$ is rigid if 
$\Hom^1(S,S)=0$. 

\begin{prop}[\cite{KQ}, Proposition 5.4]
\label{simple_objects}
Let $\D$ be a triangulated category and $\H \subset \D$ be a finite heart. 
For any rigid simple object $S \in \H$, the sets of simple objects in the simple tilted hearts 
$\H_S^{\sharp}$ and $\H_S^{\flat}$
are given by 
\begin{align*}
\Sim \H_{S}^{\sharp} &= \{S[1]\} \cup \{\,\psi^{\sharp}_S(X)\,\vert\,X \in \Sim \H,\,X \neq S\,\} \\
\Sim \H_{S}^{\flat} &= \{S[-1]\} \cup \{\,\psi^{\flat}_S(X)\,\vert\,X \in \Sim \H,\,X \neq S\,\} \\
\end{align*}
where
\begin{align*}
\psi^{\sharp}_S(X) &= \mathrm{Cone}(X \to S[1] \otimes \Hom^1(X,S)^*)[-1] \\
\psi^{\flat}_S(X) &= \mathrm{Cone}(S[-1] \otimes \Hom^1(S,X) \to X). 
\end{align*}
In addition, $\H_S^{\sharp}$ and $\H_S^{\flat}$ are also finite.
\end{prop}

\subsection{Exchange graphs}
In this section, we summarize 
basic properties of finite hearts in $\DQ$ showed in \cite{KQ}. 
\begin{defi}
The oriented graph with labeled edges 
$\EG(\DQ)$, called the exchange graph of 
$\DQ$, is defined as follows. 
The set of vertices in $\EG(\DQ)$ consists of all finite 
hearts of $\DQ$, and 
the set of labeled edges consists of all forward tiltings, i.e. 
two hearts $\H$ and $\H^{\prime}$ are connected by the arrow $\H \xto{S} \H^{\prime}$ if there is some rigid simple object $S \in \H$ such that  $\H_S^{\sharp}=\H^{\prime}$.
The connected subgraph $\EG^{\circ}(\DQ) \subset \EG(\DQ)$ is 
defined to be  
the connected component which contains the standard heart $\H_0$. 
\end{defi}

For $\H \in \EG^{\circ}(\DQ)$ and 
a rigid simple object $S \in \H$, we define inductively 
\begin{equation*}
\H^{m \sharp}_S:=(\H^{(m-1)\sharp}_S)_{S[m-1]}^{\sharp}, 
\quad \H^{m \flat}_S:=(\H^{(m-1)\flat}_S)_{S[-m+1]}^{\flat}
\end{equation*}
for $m \ge 1$.

\begin{defi}
We say that a heart $\H$ is 
\begin{itemize}
\item
spherical, if all its simples are spherical, 
\item
monochromatic, if for any simples $S \neq T \in \H$, 
$\Hom^{\bullet}(S,T)$ is concentrated in a single positive degree.  
\end{itemize}
\end{defi}

\begin{prop}[\cite{KQ}, Corollary 8.5]
\label{heart_property}
Every heart in $\EG^{\circ}(\DQ)$ is finite, spherical and 
monochromatic. Let  
$\H \in \EG^{\circ}(\DQ)$ a finite heart and 
consider simple objects 
$\Sim \H=\{T_1,\dots,T_n\} $. 
Then spherical twists $\Phi_{T_1},\dots,\Phi_{T_n}$ generate the 
Seidel-Thomas braid group $\Sph(\DQ)$. In addition, 
we have
\begin{equation*}
\H^{(N-1)\sharp}_{T_i}=\Phi^{-1}_{T_i}(\H), \quad \H^{(N-1)\flat}_{T_i}=\Phi_{T_i}(\H).
\end{equation*}
\end{prop}

The result implies that 
the action of $\Sph(\DQ)$ on $\EG(\DQ)$ preserves 
the connected component $\EG^{\circ}(\DQ)$.  
By using the result in Section \ref{Koszul}, 
we can show that this action is free.
\begin{prop}
\label{free_action}
The action of $\Br(\DQ)$ on $\EG^{\circ}(\DQ)$ is free.
\end{prop}
{\bf Proof.}
If $\Phi \in \Br(\DQ)$ acts trivially on 
$\EG^{\circ}(\DQ)$, 
then $\Phi$ preserves the vertex $\H_0$ and 
labeled edges connected to $\H_0$. This implies 
that $\Phi$ preserves simple objects $R_1,\dots,R_n$ in $\H_0$. Hence by 
Corollary \ref{tensor_id2}, it follows that $\Phi \cong \id$.
\endproof

Let $\H_i \,(i=1,2)$ be the hearts of bounded t-structures 
$\sli_i\,(i=1,2)$. 
We denote by
\begin{equation*}
\H_1 \le \H_2
\end{equation*}
if they satisfy the condition $\sli_2 \subset \sli_1$. 

We define the full subgraph of $\EG^{\circ}(\DQ)$ by
\begin{equation*}
\EG_N^{\circ}(\H_0):=\{\,\H \in \EG^{\circ}(\DQ)\,\vert\,
\H_0 \le \H \le \H_0 [N-2] \,\}.
\end{equation*}

Let $\H \in \EG(\DQ)$ be a finite heart and  
$S \in \H$ be a simple object. 
Since  
\begin{equation*}
\Phi(\H^{\sharp}_S)= (\Phi(\H))^{\sharp}_{\Phi(S)}
\end{equation*}
for any autoequivalence $\Phi \in \Aut(\DQ)$, 
the element of $\Aut(\DQ)$ acts on $\EG(\DQ)$ as an automorphism 
of the oriented graph. 
Denote by $\Aut^{\circ}(\DQ)$ the subgroup consisting of autoequivalences 
which preserve the connected component 
$\EG^{\circ}(\DQ)$.

The following implies that 
$\EG_N^{\circ}(\H_0)$ is a 
fundamental domain of 
the action of $\Sph(\DQ)$ on $\EG^{\circ}(\DQ)$.

\begin{thm}[\cite{KQ}, Theorem 8.6]
\label{fundamental}
There is a map of oriented graphs 
\begin{equation*}
p_0 \colon \EG_N^{\circ}(\H_0) \lto \EG^{\circ}(\DQ)  \slash \Sph(\DQ)
\end{equation*}
such that $p_0$ is embedding and bijective on vertices. 
\end{thm}

\subsection{Cluster categories and cluster tilting objects}
In this section, we introduce the notion of an $(N-2)$-cluster 
category, an $(N-2)$-cluster tilting object and its mutation. 
In particular, $(N-2)$-cluster tilting objects will play the role of intermediate between $N$-angulations 
and hearts of the derived category of the Ginzburg $\CY_N$ dga. 
For details of cluster categories, we refer to \cite[Section 4]{KQ}.

Let $Q$ be an acyclic quiver and  $D^b(kQ)$ 
be the bounded derived category 
of finite dimensional modules over the path algebra $kQ$. 
We denote by $\tau$ the Auslander-Reiten translation of $D^b(kQ)$. 
The $(N-2)$-cluster category $\mathcal{C}_{N-2}(Q)$ is 
defined to be the orbit category 
\begin{equation*}
\mathcal{C}_{N-2}(Q):=D^b(kQ) \slash \left<F^m \vert m \in \Z \right>
\end{equation*}
where $F:=\tau^{-1} \circ [N-2] \in \Aut(D^b(kQ))$. 
Assume that $Q$ has $n$ vertices. 
An object $P \in \mathcal{C}_{N-2}(Q)$ is called 
the $(N-2)$-cluster tilting object 
if 
\begin{itemize}
\item $\Hom_{\mathcal{C}_{N-2}(Q)}(P,P[i])=0$ 
for $i=1,\dots,N-2$,
\item 
$P$ has $n$ indecomposable 
direct summands, up to isomorphism.
\end{itemize}
Assume that an $(N-2)$-cluster tilting object 
$P$ decomposes as $P= \oplus_{i=1}^n P_i$. 
Then, one can construct a new 
$(N-2)$-cluster tilting object $\mu_i^{\sharp} P$ 
(respectively $\mu_i^{\flat} P$), called 
the forward (respectively backward) mutation of $P$ 
at the $i$-th object $P_i$.
For the construction of $\mu_i^{\sharp} P$ and 
$\mu_i^{\flat} P$, we refer to \cite[Section 4]{KQ}.

\begin{defi}
The cluster exchange graph $\mathrm{CEG}_{N-2}(Q)$ is 
an oriented graph whose vertices are $(N-2)$-cluster 
tilting objects and whose labeled edges are forward mutations, i.e. 
two tilting objects $P$ and $P^{\prime}$ are connected by 
the arrow $P \xto{P_v} P^{\prime}$ if $P^{\prime}=\mu_v^{\sharp}P$ 
where $P_v$ is an indecomposable direct summand 
of $P=\oplus_{i=1}^n P_i$. 
\end{defi}

\begin{thm}[\cite{KQ}, Theorem 8.6]
\label{KQ_correspondence}
There is an isomorphism of oriented graphs with labeled edges
\begin{equation*}
 \EG^{\circ}(\DQ) \slash \Br(\DQ)
\cong \mathrm{CEG}_{N-2}(Q).
\end{equation*}
\end{thm}

\begin{rem}
Let $\H \in \EG^{\circ}(\DQ)$ be a finite heart and $P \in \mathrm{CEG}_{N-2}(Q)$ be 
the corresponding $(N-2)$-cluster tilting object given by the above isomorphism. 
Since both $\EG^{\circ}(\DQ) \slash \Br(\DQ)$
and  $\mathrm{CEG}_{N-2}(Q)$ are $n$-regular graphs,
  the above isomorphism 
also gives the bijection between all simple objects $S_1,\dots,S_n$ in $\H$ and 
all indecomposable direct summands $P_1,\dots,P_n$ of $P$:
\begin{equation*}  
\{S_1,\dots,S_n\} \cong \{P_1,\dots,P_n\},
\quad S_i \mapsto P_i.
\end{equation*}
\end{rem}

\subsection{Colored quivers and mutation}
\label{colored_quiver_lattice}
Following 
\cite[Section 2]{BT}, we introduce the notion of a colored quiver and its mutation. 

An $(N-2)$-colored quiver $\CQ$ consists of vertices $\{1,\dots,n\}$ and 
colored arrows $i \xto{(c)}j $, where $c \in \{0,1,\dots,N-2\}$. 
Write by $q_{ij}^{(c)}$ the number of arrows from $i$ to $j$ of color $(c)$. 

In addition, we assume the following conditions:
\begin{itemize}
\item[(1)] No loops, $q^{(c)}_{ii}=0$ for all $c$.
\item[(2)] Monochromaticity, if $q^{(c)}_{ij}\neq 0$, then $q^{(c^{\prime})}=0$ for $c \neq c^{\prime}$.
\item[(3)] Skew-symmetry, $q_{ij}^{(c)}=q^{(N-2-c)}_{ji}$.
\end{itemize}

Let $\CQ$ be a $(N-2)$-colored quiver satisfying additional conditions 
and fix a vertex $v$ in $\CQ$. 
Then a new $(N-2)$-colored quiver 
$\mu^{\sharp}_v \CQ=(\tilde{q}^{(c)}_{ij})$ is defined by 
\begin{equation*}
\tilde{q}^{(c)}_{ij}= 
\begin{cases}
q_{ij}^{(c-1)}   & \text{if } v=i \\
q_{ij}^{(c+1)}   & \text{if } v=j \\
\max\{\,0,q_{ij}^{(c)}  - \sum_{t \neq c}q_{ij}^{(t)}
 +( q_{iv}^{(c)}-q_{iv}^{(c+1)} )q_{vj}^{(0)} \\
 \,\,\,\quad\qquad\quad\qquad\qquad\qquad +q_{iv}^{(N-2)}( q_{vj}^{(c)}-q_{vj}^{(c-1)} )
\} &\text{otherwise}.
\end{cases}
\end{equation*}
We call this operation 
the forward colored quiver mutation of 
$\CQ$ at a vertex $v$. The backward mutated colored quiver 
$\mu^{\flat}_v \CQ$ can be similarly defined. 
(Our notation is different to the notation in \cite{BT}. 
The backward mutation $\mu^{\flat}_v \CQ$ corresponds to  
the colored quiver mutation $\mu_v \CQ$ 
in \cite[Section 2]{BT}.)


For each $(N-2)$-cluster tilting object $P$, we can associate 
an $(N-2)$-colored quiver $\CQ(P)$ 
whose vertices $\{1,\dots,n\}$ correspond to 
indecomposable direct summands of $P=\oplus_{i=1}^n P_i$. 
The definition of $\CQ(P)$ can be found 
in \cite[Section 2]{BT} or 
\cite[Section 4]{KQ}, though the precise definition will not be used later.

\begin{thm}[\cite{BT}, Theorem 2.1]
\label{BT}
Let $P$ be an $(N-2)$-cluster tilting object and $v$ be 
a vertex of $\CQ(P)$ corresponding to an indecomposable 
direct summand $P_v$ of $P$. Then we have 
\begin{equation*}
\mu^{\sharp}_v \CQ(P) =\CQ(\mu^{\sharp}_v P),\quad
\mu^{\flat}_v \CQ(P) =\CQ(\mu^{\flat}_v P).
\end{equation*}
\end{thm}

Next we consider colored quivers associated to  
finite hearts in $\EG^{\circ}(\DQ)$.
\begin{defi}
\label{def_Hquiver}
For a finite heart $\H \in \EG^{\circ}(\DQ)$, 
we define the colored quiver $\CQ(\H)$ whose 
vertices $1,\dots,n$ correspond to 
simple objects $S_1,\dots,S_n \in \Sim \H$ and 
the number of arrows $q_{ij}^{(c)}$ from $i$ to $j$ of color $(c)$ 
is given by $q_{ij}^{(c)}:=\dhom^{c+1}(S_i,S_j)$ for $i \neq j$ 
and $q_{ii}^{(c)}:=0$. 
\end{defi}
By definition, $\CQ(\H)$ has no loops. 
Further, the $\CY_N$ property of $\DQ$ implies that 
$\CQ(\H)$ is skew-symmetric 
and Proposition \ref{heart_property} 
implies the monochromaticity of $\CQ(\H)$. 
Therefore, the colored quiver $\CQ(\H)$ satisfies 
additional conditions. 

\begin{rem}
For a finite heart $\H$, we can associate a graded quiver 
$\CQ^{\epsilon}(\H)$, 
called the Ext quiver of  $\H$ (Definition 6.2 in \cite{KQ}), 
whose vertices 
are simple objects of $\H$ and whose degree $k$ arrows 
from $S_i$ to $S_j$ correspond to a basis of 
$\Hom_{\DQ}(S_i,S_j[k])$. 
On the other hand, if we consider the augmented graded quiver $\CQ^+(\H)$ 
(Definition 6.1 in \cite{KQ})
of the colored quiver $\CQ(\H)$, then we have 
$\CQ^{\epsilon}(\H)=\CQ^+(\H)$. 
\end{rem}

The following result is stated in terms of 
the Ext quiver and the augmented quiver in 
the original paper \cite{KQ}, however 
here we state in terms of the colored quiver.  
\begin{thm}[\cite{KQ}, Theorem 8.10]
\label{KQ_quiver}
Let $\H \in \EG^{\circ}(\DQ)$ be a finite heart 
and $P \in \mathrm{CEG}_{N-2}(Q).$ be 
the corresponding $(N-2)$-cluster tilting object given by 
Theorem \ref{KQ_correspondence}. 
Then we have 
\begin{equation*}
\CQ(\H)=\CQ(P).
\end{equation*}
\end{thm}

By Theorem \ref{KQ_correspondence}, Theorem \ref{BT} and Theorem \ref{KQ_quiver}, 
we have the following corollary.
\begin{cor}
\label{mutation}
Let $\H \in \EG^{\circ}(\DQ)$ be a finite heart 
and $\CQ(\H)$ be 
the corresponding colored quiver. 
Then for any simple object $S_v \in \H$, we have
\begin{equation*}
\mu^{\sharp}_v\CQ(\H)= \CQ(\H_{S_v}^{\sharp}), \quad
\mu^{\flat}_v\CQ(\H)= \CQ(\H_{S_v}^{\flat}).
\end{equation*}
\end{cor}
{\bf Proof.} 
Let $P$ be the $(N-2)$-cluster tilting object corresponding to $\H$ 
given by Theorem \ref{KQ_correspondence}. 
Again by Theorem \ref{KQ_correspondence}, the $(N-2)$-cluster tilting object 
corresponding to $\H_{S_v}^{\sharp}$ is the forward mutation $\mu_v^{\sharp}P$. 
Then, Theorem \ref{BT} and Theorem \ref{KQ_quiver} imply that
\begin{equation*}
\mu^{\sharp}_v\CQ(\H)  
=\mu^{\sharp}_v\CQ(P) 
=\CQ(\mu_v^{\sharp}P) 
=\CQ(\H_{S_v}^{\sharp}).
\end{equation*}
\endproof

\subsection{Colored quiver lattices}
\label{sec_mutation_K}
In this section, we define a lattice associated to a 
colored quiver. 

For a $(N-2)$-colored quiver  $\CQ$ with vertices $\{1,\dots,n\}$, 
we introduce a free abelian group $L_{\mathcal{Q}}$ 
generated by vertices of $\mathcal{Q}$. 
Denote by $\alp_1,\dots,\alp_n$ generators corresponding to vertices $1,\dots,n$. 
Then we have 
\begin{equation*}
L_{\mathcal{Q}} = \bigoplus_{i=1}^n \Z \,\alp_i.
\end{equation*} 
Define the bilinear form $\left<\,,\,\right> \colon \LQ \times \LQ \lto \Z$ by 
\begin{equation*}
\left<\alp_i,\alp_j \right> := \delta_{ij} + (-1)^{N}\delta_{ij} -\sum_{c=0}^{N-2} (-1)^c q_{ij}^{(c)}.
\end{equation*}
Since $q_{ij}^{(c)}= q_{ji}^{(N-2-c)}$, the bilinear form 
$\left<\,,\,\right>$ is symmetric if $N$ is even, 
and skew-symmetric if $N$ is odd.

Let $\CQ(\H)$ be the colored quiver 
of a heart $\H \in \EG^{\circ}(\DQ)$ with 
simple objects $S_1,\dots,S_n$.
We write by $\alp_i$ the basis of 
colored quiver lattice $L_{\CQ(\H)}$ 
which corresponds to a vertex $i$ of $\CQ(\H)$. 
Then we have the following obvious result. 

\begin{lem}
\label{quiver_K}
There is an isomorphism of abelian groups
\begin{equation*}
\lam_{\H} \colon L_{\CQ(\H)} \lto K(\H) \cong K(\D_Q^N)
\end{equation*}
defined by $\alp_i \mapsto [S_i]$. The map $\lam$ 
takes the bilinear form $\left<\,,\,\right>$ on $L_{\CQ(\H)}$ 
to the Euler form $\chi$ on $K(\D_Q^N)$.
\end{lem}
{\bf Proof.}
Since 
\begin{equation*}
 K(\H) \cong \bigoplus_{i=1}^n \Z [S_i],
\end{equation*}
$\lam_{\H}$ is an isomorphism. 
Recall from Definition \ref{def_Hquiver} 
that $q_{ij}^{(c)}=
\hom^{c+1}(S_i,S_j)$. In addition, since $\H$ is 
spherical, $\hom^0(S_i,S_j)=\hom^N(S_i,S_j)=\delta_{ij}$. 
Hence the result follows.
\endproof


Let $\H^{\sharp}_{S_v}$ be the forward simple tilted heart  
and $S_1^{\prime},\dots,S_1^{\prime}$ be 
its simple objects.
By Proposition \ref{simple_objects}, we have 
\begin{equation}
\label{eq_simple}
S_j^{\prime}= 
\begin{cases}
S_v[1] & \text{if } j=v \\
\psi^{\sharp}_{S_v}(S_j)   & \text{if } j \neq v 
\end{cases}
\end{equation}
where 
\begin{equation*}
\psi^{\sharp}_{S_v}(S_j) = \mathrm{Cone}(S_j \to S_v[1] \otimes \Hom^1(S_v,S_j)^*)[-1].
\end{equation*}

Let $\mu_v^{\sharp}\CQ(\H)$ be the forward colored 
quiver mutation of $\CQ(\H)$ at $v$.

Consider two colored quiver lattices
\begin{equation*}
L_{\CQ(\H)}=\bigoplus_{j=1}^n \Z\, \alp_j,\quad 
L_{\mu^{\sharp}_v \CQ(\H)}=\bigoplus_{j=1}^n \Z\, \alp_j^{\prime}.
\end{equation*}

We define a linear map 
$F_v \colon L_{\mu^{\sharp}_v \CQ(\H)} \to 
L_{\CQ(\H)}$ by 
\begin{equation*}
F_v(\alp_j^{\prime}):=
\begin{cases}
\,-\alp_v &\text{if}\quad j=v \\
\,\alp_j + q_{vj}^{(0)}\alp_v 
&\text{if}\quad j \neq v.
\end{cases}
\end{equation*}
Since $\mu_v^{\sharp}\CQ(\H)=\CQ(\H_{S_v}^{\sharp})$ 
by Corollary \ref{mutation}, we have 
$L_{\mu_v^{\sharp}\CQ(\H)} \cong L_{\CQ(\H_{S_v}^{\sharp})}$. 
Hence from Lemma \ref{quiver_K}, we have an isomorphism 
\begin{equation*}
\lam_{\H_{S_v}^{\sharp}} \colon L_{\mu_v^{\sharp}\CQ(\H)} 
\lto K(\H_{S_v}^{\sharp}) \cong K(\DQ)
\end{equation*}
defined by $\alp_j^{\prime} \mapsto [S^{\prime}_j]$. 

With the above notations, we have the following diagram.

\begin{lem}
\label{diagram1}
Let $\lam_{\H} \colon L_{\CQ(\H)} \to K(\D_n^N)$ and 
$\lam_{\H_{S_v}^{\sharp}} \colon L_{\mu_v^{\sharp} \CQ(\H)} \to K(\D_n^N)$ be the isomorphisms 
given by Lemma \ref{quiver_K}. Then the diagram
\begin{equation*}
\xymatrix{ 
L_{\mu^{\sharp}_v \CQ(\H)} \ar[rd]_{\lam_{\H_{S_v}^{\sharp}}} 
\ar[rr]^{F_v} &&
L_{\CQ(\H)} \ar[ld]^{\lam_{\H}} \\ 
& K(\DQ) &
}
\end{equation*}
commutes. 
\end{lem}
{\bf Proof.}
First note that by the equation (\ref{eq_simple}), we have 
\begin{equation*}
[S_j^{\prime}]= 
\begin{cases}
-[S_v] & \text{if } j=v \\
\,[S_j]+ \dhom^1(S_v,S_j)[S_v]  
& \text{if } j \neq v 
\end{cases}
\end{equation*}
in $K(\DQ)$. Recall from Definition \ref{def_Hquiver} 
that $q_{vj}^{(0)}=\dhom^1(S_v,S_j)$. 
Since $\lam_{\H}(\alp_j)=[S_j]$ and 
$\lam_{\H_{S_v}^{\sharp}}(\alp_j^{\prime})=[S_j^{\prime}]$, 
we have the result. 
\endproof

\section{Bridgeland stability conditions}
\label{sec:stab}
\subsection{Spaces of stability conditions}

In this section, we recall the notion of 
stability conditions on triangulated categories 
introduced by T. Bridgeland \cite{Br1} and collect 
some results for spaces of stability conditions. 

\begin{defi}
Let $\D$ be a triangulated category and $K(\D)$ be its $K$-group. 
A stability condition $\sigma = (Z, \sli)$ on $\D$ consists of 
a group homomorphism $Z \colon K(\D) \to \C$ called central charge and 
a family of full additive subcategories $\sli (\phi) \subset \D$ for $\phi \in \R$ 
satisfying the following conditions:
\begin{itemize}
\item[(a)]
if  $0 \neq E \in \sli(\phi)$, 
then $Z(E) = m(E) \exp(i \pi \phi)$ for some $m(E) \in \R_{>0}$,
\item[(b)]
for all $\phi \in \R$, $\sli(\phi + 1) = \sli(\phi)[1]$, 
\item[(c)]if $\phi_1 > \phi_2$ and $A_i \in \sli(\phi_i)\,(i =1,2)$, 
then $\Hom_{\D}(A_1,A_2) = 0$,
\item[(d)]for $0 \neq E \in \D$, there is a finite sequence of real numbers 
\begin{equation*}
\phi_1 > \phi_2 > \cdots > \phi_m
\end{equation*}
and a collection of exact triangles
\begin{equation*}
0 =
\xymatrix @C=5mm{ 
 E_0 \ar[rr]   &&  E_1 \ar[dl] \ar[rr] && E_2 \ar[dl] 
 \ar[r] & \dots  \ar[r] & E_{m-1} \ar[rr] && E_m \ar[dl] \\
& A_1 \ar@{-->}[ul] && A_2 \ar@{-->}[ul] &&&& A_m \ar@{-->}[ul] 
}
= E
\end{equation*}
with $A_i \in \sli(\phi_i)$ for all $i$.
\end{itemize}
\end{defi}

It follows from the definition that the subcategory 
$\sli(\phi) \subset \D$ is an abelian category 
(see Lemma 5.2 in \cite{Br1}). 
Nonzero objects in $\sli(\phi)$ are called semistable of phase $\phi$ in $\sigma$ and simple objects 
in $\sli(\phi)$ are called stable of phase $\phi$ in $\sigma$. 

For a stability condition $\sigma = (Z,\sli)$, 
we introduce the set of semistable classes $\class^{\mathrm{ss}}(\sigma) \subset K(\D)$ by
\begin{equation*}
\class^{\mathrm{ss}}(\sigma) :=\{\,\alp \in 
K(\D)\,\vert\,\text{there exists a semistable object } E \in \D \text{ in $\sigma$ such that } [E] = \alp\,\}.
\end{equation*} 

We always assume that our stability conditions satisfy the additional assumption,  
called the support property in \cite{KS}. 
\begin{defi}
\label{support}
Let $\norm{\,\cdot\,}$ be some norm on $K(\D) \otimes \R$. A stability condition $\sigma=(Z,\sli)$ 
satisfies the support property if there is a some constant $C >0$ such that 
\begin{equation*}
C \cdot \abv{Z(\alp)} > \norm{\alp}
\end{equation*}
for all $\alp \in \class^{\mathrm{ss}}(\sigma)$. 
\end{defi}

Let $\Stab(\D)$ be the set of all stability conditions on $\D$ with the support property. 

In \cite{Br1}, Bridgeland introduced a natural topology on the space $\Stab(\D)$ 
characterized by the following theorem. 
\begin{thm}[\cite{Br1}, Theorem 1.2]
\label{localiso}
The space $\Stab(\D)$ has the structure of a complex manifold and 
the projection map of central charges
\begin{equation*}
\mathcal{Z} \colon \Stab(\D) \longrightarrow \Hom_{\Z}(K(\D),\C)
\end{equation*}
defined by $(Z,\sli) \mapsto Z$ is a local isomorphism of complex manifolds 
onto an open subset of $\Hom_{\Z}(K(\D),\C)$. 
\end{thm} 

In \cite{Br1}, Bridgeland gave the alternative description of a stability condition on $\D$ 
as a pair of a bounded t-structure and a central charge on its heart. 

\begin{defi}
Let $\H$ be an abelian category and let $K(\H)$ 
be its Grothendieck group. 
A central charge on $\H$ is 
a group homomorphism $Z \colon K(\H) \to \mathbb{C}$ such that for any nonzero object
$0 \neq E \in \H$, the complex number $Z(E)$ lies in semi-closed upper half-plane 
$H = \{\, r e^{i \pi \phi}\ \in \mathbb{C} \,
 \vert \, r \in \mathbb{R}_{>0}, \phi \in (0,1] \,\}$.
\end{defi}

The real number 
\begin{equation*}
\phi(E):= \frac{1}{\pi}\,\arg Z(E) \in (0,1]
\end{equation*}
for $0 \neq E \in \H$ is called the phase of $E$. 

A nonzero object $0 \neq E \in \H$ is said to be $Z$-(semi)stable if 
for any nonzero proper subobject $0 \neq A \subsetneq E$ satisfies $\phi(A)  < (\le) \phi(E)$. 

For a central charge $Z \colon K(\H) \to \C$, the set of semistable classes can be defined similarly. 

\begin{prop}[\cite{Br1}, Proposition 5.3]
\label{abel_stability}
Let $\D$ be a triangulated category. To give a stability condition on $\D$ is equivalent 
to giving the heart $\H \subset \D$ 
of a bounded structure on $\D$ and a central charge with the Harder-Narasimhan property on $\H$. 
\end{prop}
For details of the Harder-Narasimhan property 
(HN property), we refer 
to \cite[Section 2]{Br1}. 

We denote by $\Stab(\H)$ the set of central charges on 
the heart  $\H \subset \D$  
with the HN property and the support property.

\subsection{Simple tilting and gluing}
Here by using Proposition \ref{abel_stability},  
we construct stability conditions on a finite heart.

Let $\H \subset \D$ be a finite heart with simple objects $\{S_1,\dots,S_n\}$. 
Then $K(\H) \cong \oplus_{i=1}^n \Z [S_i]$. 
For any point $(z_1,\dots,z_n) \in H^n$, the central charge $Z \colon K(\H) \to \C$ is 
constructed by $Z(S_i) :=z_i$. Conversely, for a given central charge 
$Z \colon K(\H) \to \C$, the complex number 
$Z(S_i)$ lies in $H$ for all $i$. Hence $Z$ determines 
a point $(z_1,\dots,z_n) \in H^n$ where $z_i :=Z(S_i)$.
As a result, the set of central charges on $\H$ is 
isomorphic to $H^n$.

\begin{lem}[\cite{Br4}, Lemma 5.2]
Let $Z \colon K(\H) \to \C$ be a central charge given by the above construction. 
Then $Z$ has a Harder-Narasimhan property. 
In particular, we have 
\begin{equation*}
\Stab(\H) \cong H^n. 
\end{equation*}
\end{lem}

Let $\H \subset \D$ be a finite rigid heart with simple objects $\{S_1,\dots,S_n\}$.
From the above results, there is a natural inclusion 
\begin{equation*}
\Stab(\H) \subset \Stab(\D).
\end{equation*} 
From such a heart, we can construct the another hearts $\H_{S_i}^{\sharp}$(or $\H_{S_i}^{\flat}$) 
by simple tilting as in Section \ref{sec_tilting}. 
The relationship between $\Stab(\H)$ and 
$\Stab(\H_{S_i}^{\sharp})$(or $\Stab(\H_{S_i}^{\flat})$) in $\Stab(\D)$ is given by the following lemma.

\begin{lem}[\cite{BrSm}, Lemma 7.9]
\label{gluing}
Let $\sigma=(Z,\sli) \in \Stab(\D)$ lies on a unique codimension one boundary of the 
region $\Stab(\H)$ so that $\Im Z(S_i)=0$ for a unique simple object $S_i$. 
Then there is a neighborhood $\sigma \in U \subset \Stab(\D)$ such that one 
of the following holds:
\begin{itemize}
\item[(1)] $Z(S_i) \in \R_{>0}$ and $U \subset \Stab(\H) \cup \Stab(\H_{S_i}^{\sharp})$,
\item[(2)] $Z(S_i) \in \R_{<0}$ and $U \subset \Stab(\H) \cup \Stab(\H_{S_i}^{\flat})$.
\end{itemize}
\end{lem}

Consider $\D=\D_n^N$, which is the derived category of 
finite total dimensional dg modules over the 
Ginzburg dg algebra $\Gamma_N \QA$ associated to 
the $A_n$-quiver (see section \ref{Ginzburg}). 
By Proposition \ref{standard_heart}, there is a special 
heart $\H_0 \subset \D_n^N$, called the standard heart. 
Since $\H_0$ is finite, we have $\Stab(\H_0) \cong H^n$. 
\begin{defi}
The distinguished connected component 
$\Stab^{\circ}(\D_n^N) \subset \Stab(\D_n^N)$ 
is defined to be  
the connected component which contains the 
connected subset $\Stab(\H_0)$.  
\end{defi}

By using the result in \cite{Wol}, Qiu described the 
space $\Stab^{\circ}(\D_n^N)$ as follows. 

\begin{prop}[\cite{Q}, Corollary 5.3]
\label{prop_gluing}
The distinguished connected component $\Stab^{\circ}(\D_n^N)$ 
is given by 
\begin{equation*}
\Stab^{\circ}(\D_n^N) =
\bigcup_{\H \in \EG^{\circ}(\D_n^N) }\Stab(\H).
\end{equation*}
\end{prop}

\subsection{Group actions}
\label{sec_group}
For the space $\Stab(\D)$, we introduce two commuting group actions. 

First we consider the action of $\Aut(\D)$ on $\Stab(\D)$. 
Let $\Phi \in \Aut(\D)$ be an autoequivalence of $\D$ and 
 $(Z,\sli) \in \Stab(\D)$ be a stability condition on $\D$. Then the element 
$\Phi \cdot (Z,\sli) =(Z^{\prime},\sli^{\prime})$ is defined by  
\begin{equation*}
Z^{\prime}(E):=Z(\Phi^{-1}(E)), \quad \sli^{\prime}(\phi) := \Phi(\sli(\phi)),
\end{equation*}
where $E \in \D$ and $\phi \in \R$.
 
Second we define the $\C$-action on $\Stab(\D)$. 
For $t \in \C$ and $(Z,\sli) \in \Stab(\D)$, 
the element $t \cdot (Z,\sli) =(Z^{\prime},\sli^{\prime})$ is defined by 
\begin{equation*}
Z^{\prime}(E):= e^{-i \pi t} \cdot Z(E), \quad \sli^{\prime}(\phi) :=\sli(\phi + \Re (t))
\end{equation*}
where $E \in \D$ and $\phi \in \R$. 
Clearly, this action is free. 

Consider the category $\D=\D_n^N$. 
Recall from Section \ref{sec_spherical} that there is a 
subgroup $\Sph(\D_n^N) \subset \Aut(\D_n^N)$ consisting  
of spherical twists. 
By Proposition \ref{free_action} and 
Proposition \ref{prop_gluing}, 
we have the following. 
\begin{lem}
\label{free}
The action of $\Sph(\D_n^N)$ on $\Stab^{\circ}(\D_n^N)$ 
is free and properly discontinuous. 
\end{lem}
{\bf Proof.} 
Proposition \ref{free_action} and 
Proposition \ref{prop_gluing} directly 
imply that the action of 
$\Sph(\D_n^N)$ on $\Stab^{\circ}(\D_n^N)$ is free. 

We show that the action is properly discontinuous. 
Again by Proposition \ref{prop_gluing}, it is sufficient to prove 
that for a stability condition $\sigma \in \Stab(\H) $ with  
$\H \in \EG^{\circ}(\D_n^N)$, 
there is an open neighborhood $\sigma \in U \subset 
\Stab^{\circ}(\D_n^N)$ such that $U \cap  
\Phi( U) = \emptyset $ for $\id \ncong \Phi 
\in \Sph(\D_n^N)$. 
Let $\sigma \in \Int(\Stab(\H) )$ be a stability condition 
in the interior of $\Stab(\H) $ and 
take an open neighborhood 
$\sigma \in U \subset \Int(\Stab(\H) ) $. 
Since $\H \neq \Phi(\H)$ for $\id \ncong \Phi 
\in \Sph(\D_n^N)$, 
$\Stab(\H) \cap \Stab(\Phi(\H))=\emptyset$. Hence 
$U \cap \Phi(U)=\emptyset$. For a general element 
$\sigma \in \Stab(\H) $, there is some complex number 
$t \in \C$ such that $t \cdot \sigma \in \Int(\Stab(\H))$. 
Since the action of $\C$ commutes with the action of 
$\Br(\DQ)$, we can apply the same argument as for $\sigma$ 
an element in the interior of $\Stab(\H)$. 
\endproof

\section{$N$-angulations and hearts}
\label{sec_6}
\subsection{$N$-angulations of polygons}
\label{$N$-angulation}
In \cite{BaMa}, they give a geometric description of the  
$(N-2)$-cluster category of type $A_n$.  
By using their model, we give 
a geometric realization of the cluster exchange graph 
of an $(N-2)$-cluster category of type $A_n$. 
 
Let $\Pi_{\d}$ be a $\d$-gon.  
An $(N-2)$-diagonal of $\Pi_{\d}$ is a diagonal which divides $\Pi_{\d}$ into an $((N-2)i+2)$-gon 
and an $((N-2)(n+1-i)+2)$-gon ($i=1,\dots,n$). 
A collection of $(N-2)$-diagonals is called non-crossing if they 
intersect only at the boundary of $\Pi_{\d}$. 
\begin{defi}
An $N$-angulation $\Delta$ of $\Pi_{\d}$ is the maximal set of 
non-crossing $(N-2)$-diagonals of $\Pi_{\d}$. 
\end{defi}
Note that the number of $(N-2)$-diagonals in the 
$N$-angulation $\Delta$ is just $n$ 
and they divide $\Pi_{\d}$ into $(n+1)$  $N$-gons. 

For any $(N-2)$-diagonal $\delta \in \Delta$, 
there are just two $N$-gons whose common edge is $\delta$. 
By removing $\delta$, we have a $(2N-2)$-gon. We call $\delta$ a 
diameter of the $(2N-2)$-gon.  

For an $N$-angulation $\Delta$ and 
an $(N-2)$-diagonal $\delta \in \Delta$, 
we introduce the flip operation to make new $N$-angulations $\mu_{\delta}^{\sharp}\Delta$ (called the forward flip 
w.r.t. $\delta$) and 
$\mu_{\delta}^{\flat}\Delta$ (called the backward flip w.r.t. $\delta$) 
by the following. 

For an $(N-2)$-diagonal $\delta$, there is a unique 
$(2N-2)$-gon which has $\delta$ as a diameter. 
By rotating two boundary points of $\delta$ one step 
clockwise  (respectively counterclockwise), we have 
a new $(N-2)$-diagonal $\delta^{\sharp}$ 
(respectively $\delta^{\flat}$). 
Thus, we obtain a new $N$-angulation 
\begin{equation*}
\mu_{\delta}^{\sharp}\Delta:=(\Delta \bs \{\delta\}) \cup \{\delta^{\sharp}\} 
\quad (\,\mu_{\delta}^{\flat}\Delta:=(\Delta \bs \{\delta\}) \cup \{\delta^{\flat}\}\,).
\end{equation*}

\begin{figure}[!ht]
\centering
\includegraphics[scale=0.7]{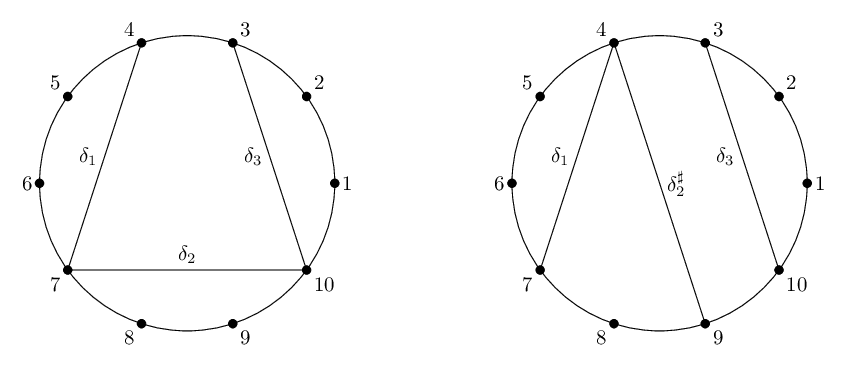} 
\caption{Example of an $N$-angulation and 
its forward flip in the case $N=4$ and $n=3$.}
\label{Rotation.fig}
\end{figure}

\begin{defi}
Define the graph $\Ang(\Pi_{\d},N)$ to be an oriented graph whose vertices are all $N$-angulations of 
$\Pi_{\d}$, and for any two $N$-angulations $\Delta$ and $\Delta^{\prime}$, the 
arrow $\Delta \xto{\delta} \Delta^{\prime}$ is given if 
there is some $(N-2)$-diagonal $\delta \in \Delta$ such that
$\Delta^{\prime}=\mu_{\delta}^{\sharp}\Delta$. 
\end{defi}


The geometric description of the $(N-2)$-cluster category of type $A_n$ 
in \cite{BaMa} implies the following correspondence. 
\begin{thm}[\cite{BaMa}, Theorem 5.6]
\label{BaMa}
Let $Q$ be an $A_n$-quiver. Then, 
there is an isomorphism of oriented graphs with labeled arrows 
\begin{equation*}
\Ang(\Pi_{\d},N) \xto{\sim} \mathrm{CEG}_{N-2}(Q).
\end{equation*}  
\end{thm}

\subsection{Correspondence between hearts and $N$-angulations}
Combining the geometric description of the $(N-2)$-cluster category of type $A_n$ (\cite{BaMa}) 
with the correspondence between the cluster exchange graph and the heart exchange graph (\cite{KQ}), 
we have the geometric description of the heart exchange graph for an $A_n$-quiver. 
From Theorem \ref{KQ_correspondence} and Theorem \ref{BaMa}, we have the following. 
\begin{cor}
\label{correspondence}
There is an isomorphism of oriented graphs with labeled arrows 
\begin{equation*}
\Ang(\Pi_{\d},N) \xto{\sim} \EG^{\circ}(\D^N_n)  \slash \Sph(\D^N_n).
\end{equation*}  
For  $\Delta \in \Ang(\Pi_{\d},N)$ and the corresponding 
heart $\H(\Delta) \in  \EG^{\circ}(\D^N_n)$ 
(which is determined up to modulo $\Sph(\D^N_n)$),
there is a canonical bijection between $(N-2)$-diagonals 
$\Delta=\{ \delta_1,\dots,\delta_n\}$ and 
simple objects $\Sim \H(\Delta)=
\{S_1,\dots,S_n\}$:
\begin{equation}
\{\delta_1,\dots,\delta_n \} \iso
\{S_1,\dots,S_n\},\quad \delta_i \mapsto S_i.
\end{equation}
In addition, we have 
\begin{equation*}
\H(\mu_{\delta_v}^{\sharp}\Delta) \equiv \H(\Delta)_{S_v}^{\sharp} \quad \mod \Sph(\D^N_n)
\end{equation*}
where $\delta_v \in \Delta$ is an $(N-2)$-diagonal and 
$S_v \in \Sim \H(\Delta)$ is the corresponding simple object. 
\end{cor} 


We fix the ambiguity of $\H(\Delta)$ as follows. 

\begin{rem}
\label{heart_fundamental}
By Theorem \ref{fundamental}, for an $N$-angulation $\Delta$ 
we can take the corresponding heart $\H(\Delta)$ uniquely 
in the fundamental domain $\EG_N^{\circ}(\H_0)$. 
In the following, we always assume that $\H(\Delta)$ is 
in $\EG_N^{\circ}(\H_0)$. 
\end{rem}

\begin{lem}
\label{lem_spherical}
Let $\Delta \in \Ang(\Pi_{\d},N)$ be an $N$-angulation and $\H(\Delta) \in \EG_N^{\circ}(\H_0)$ 
be the corresponding heart. 
Take a $(N-2)$-diagonal $\delta_v \in \Delta$ and consider 
the new $N$-angulation $\mu_{\delta_v}^{\sharp}\Delta$. 
Then there is a unique spherical twist $\Phi \in \Sph(\D_n^N)$ such that 
\begin{equation*}
\Phi \cdot \H(\mu_{\delta_v}^{\sharp}\Delta)=
\H(\Delta)_{S_v}^{\sharp}.
\end{equation*}
\end{lem}
{\bf Proof.} 
By Corollary \ref{correspondence}, such a spherical twist 
exists. In addition, since the action of $\Sph(\D_n^N)$ on $\EG^{\circ}(\D^N_n)$ is free, the result follows. 
\endproof

\subsection{Colored quivers from $N$-angulations}
Following \cite[Section 11]{BT}, we associate the  
colored quiver $\CQ(\Delta)$ to an $N$-angulation $\Delta$.

\begin{defi}
\label{colored_quiver}
For an $N$-angulation $\Delta \in \Ang(\Pi_{\d},N)$, 
define the $(N-2)$-colored quiver $\CQ(\Delta)$ whose vertices are $(N-2)$-diagonals in $\Delta$ 
and the colored arrow $\delta_i \xto{(c)} \delta_j$ is given if 
both $\delta_i$ and $\delta_j$ lie on the same $N$-gon in the $N$-angulated $\d$-gon and 
there are just $c$ edges forming the segment of the boundary of the $N$-gon from $\delta_i$ to $\delta_j$ in the 
counterclockwise direction. 
\end{defi}

\begin{prop}
\label{correspondence_quiver}
Let $\Delta \in \Ang(\Pi_{\d},N)$ be an $N$-angulation and 
$\H(\Delta) \in \EG^{\circ}(\D^N_n)$ be the corresponding heart. 
Then there is a canonical isomorphism of colored quivers
\begin{equation*}
\CQ(\Delta)=\CQ(\H(\Delta)). 
\end{equation*}
In addition, the correspondence between vertices of 
two colored quivers
is given by the correspondence 
between $(N-2)$-diagonals and simple objects of $\H(\Delta)$ in Corollary \ref{correspondence}. 
\end{prop}
{\bf Proof.} By Proposition 11.1 in \cite{BT}, the colored quiver $\CQ(\Delta)$ 
is equal to the colored quiver $\CQ(P(\Delta))$ associated to the cluster tilting object $P(\Delta)$
given by Theorem \ref{BaMa}. 
On the other hand, 
Theorem \ref{KQ_quiver} implies that the colored quiver 
$\CQ(P(\Delta))$ coincides with 
the colored quiver $\CQ(\H(\Delta))$. 
\endproof  

The definition of the colored quiver $\CQ(\Delta)$ is compatible with the mutation 
and the flip. 
\begin{lem}
\label{ang_mu}
Let $\delta \in \Delta$ be an $N$-angulation. Then we have 
\begin{equation*}
\mu^{\sharp}_\delta\CQ(\Delta)=\CQ(\mu^{\sharp}_\delta \Delta). 
\end{equation*}
\end{lem}

Let $\Delta=\{\delta_1,\dots,\delta_n\}$ be an $N$-angulation 
and $\mu^{\sharp}_{\delta_v} \Delta$ be the forward flip w.r.t. $\delta_v$. 
By Lemma \ref{correspondence_quiver} and Lemma \ref{ang_mu}, 
we have canonical isomorphisms of colored quiver lattices
\begin{align*}
L_{\CQ(\Delta)} &\cong L_{\CQ(\H(\Delta))} \\
L_{\mu^{\sharp}_{\delta_v}\CQ(\Delta)} &\cong 
L_{\CQ(\mu^{\sharp}_{\delta_v} \Delta)} \cong
L_{\CQ(\H(\mu^{\sharp}_{\delta_v} \Delta))}.
\end{align*}

Let $S_1,\dots, S_n$ and $S^{\sharp}_1,\dots,S^{\sharp}_n$ 
be simple objects in $\H(\Delta)$ and $\H(\mu^{\sharp}_{\delta_v} \Delta)$.

By using Lemma \ref{quiver_K}, we have isomorphisms of 
lattices (for simplicity, we write $\lam:=\lam_{\H(\Delta)}$ and $\lam^{\sharp}:=\lam_{\H(\mu^{\sharp}_{\delta_v} \Delta)} $)
\begin{align*}
&\lam \colon L_{\CQ(\Delta)} \to K(\H(\Delta))  \cong K(\D_n^N) \\
&\lam^{\sharp} \colon L_{\mu^{\sharp}_{\delta_v}\CQ(\Delta)} 
\to K(\H(\mu^{\sharp}_{\delta_v} \Delta)) \cong K(\D^N_n)
\end{align*}
given by $\lam(\alp_i):=[S_i]$ and 
$\lam^{\sharp}(\alp_i^{\prime}):=[S^{\sharp}_i]$.

Then we have the following diagram. 

\begin{lem}
\label{diagram11}
Let $F_v \colon L_{\mu^{\sharp}_v\CQ(\Delta)} \to L_{\CQ(\Delta)}$ 
be the linear map defined in Section \ref{sec_mutation_K} and 
$\Phi \in \Br(\D_n^N)$ be a unique spherical twist satisfying 
$\Phi \cdot \H(\mu_{\delta_v}^{\sharp}\Delta)=
\H(\Delta)_{S_v}^{\sharp}$ (see Lemma \ref{lem_spherical}).
Then the diagram
\begin{equation*}
\xymatrix{ 
L_{\mu^{\sharp}_v\CQ(\Delta)} \ar[d]_{\lam^{\sharp}} \ar[rr]^{F_v} && 
L_{\CQ(\Delta)} \ar[d]^{\lam} \\ 
K(\D_n^N) \ar[rr]^{[\Phi]} && K(\D_n^N)
}
\end{equation*}
commutes. 
\end{lem}
{\bf Proof.} 
As in Section \ref{sec_mutation_K}, define the isomorphism 
\begin{equation*}
\lam^{\prime}:=\lam_{\H(\Delta)_{S_v}^{\sharp}}
 \colon L_{\mu^{\sharp}_v\CQ(\Delta)} \to 
 K(\H(\Delta)_{S_v}^{\sharp}) \cong K(\D_n^N)
\end{equation*} 
by $\alp^{\prime}_i \mapsto [S_i^{\prime}]$ where 
$S_1^{\prime},\dots,S_n^{\prime}$ are simple objects in 
$\H(\Delta)_{S_v}^{\sharp}$. Then, since $[\Phi] \cdot 
[S_i^{\sharp}]=
[S_i^{\prime}]$ by Lemma \ref{lem_spherical}, 
we have a commutative diagram
\begin{equation*}
\xymatrix{ 
& L_{\mu^{\sharp}_v \CQ(\Delta)} \ar[ld]_{\lam^{\sharp}} 
\ar[rd]^{\lam^{\prime}} & \\ 
 K(\D_n^N) \ar[rr]^{[\Phi]} && K(\D_n^N).
}
\end{equation*}
Combining the above diagram with the diagram in 
Lemma \ref{diagram1}, we have the result. 
\endproof

\subsection{$N$-angulations from saddle-free differentials}
\label{sec_angulation}
Let $\phi \in \Q$ be a quadratic differential. 
Recall from Section \ref{critical} that $\phi$ 
defines $\d$ distinguished tangent directions at $\infty \in \P^1$ 
since $\phi$ has a pole of order $(\d+2)$ at $\infty$. 
Consider the real oriented blow-up at $\infty \in \P^1$ (we replace 
$\infty$ with $S^1$ which corresponds to tangent directions at 
$\infty$), then 
we have a disk $\Pi$. By adding $\d$ points on the boundary of $\Pi$, which correspond 
to $\d$ distinguished tangent directions, we have a $\d$-gon $\Pi_{\d}$. 

In \cite[Section 6]{GMN}, 
they constructed ideal triangulations of marked 
bordered surfaces from saddle-free differentials with simple zeros 
(see also \cite[Section 10.1]{BrSm}). Motivated by their works,  
we consider the following construction. 

\begin{lem}
\label{N-angulation}
Let $\phi \in \Q$ be a saddle-free differential. 
By taking one generic trajectory from each horizontal strip in the decomposition of $\P^1$ by $\phi$, 
we have an $N$-angulation $\Delta_{\phi}$ of $\Pi_{\d}$.
\end{lem}
{\bf Proof.} By Lemma \ref{saddle_number}, there are just $n$ horizontal strips in 
the decomposition,  
and generic trajectories taken from different horizontal strips are non-crossing. 
Hence, it is sufficient to prove that these  generic trajectories are $(N-2)$-diagonals. 
Assume that some generic trajectory divides 
$\Pi_{\d}$ into two disks $\Pi_1$ and $\Pi_2$, 
and $\Pi_1$ contain $i$ zeros. 
Then the restriction of $\phi$ on the interior of $\Pi_1$ is 
also saddle-free with $i$ zeros. 
Therefore, the foliation in  
$\Pi_1$ has $(i-1)$ horizontal strips and $Ni$ separating trajectories. 
By the equation (\ref{eq3.2}), the number of half planes in $\Pi_1$ is $((N-2)i +2)$, 
and this implies that $\Pi_1$ is $((N-2)i+2)$-gon. 
\endproof

\begin{figure}[!ht]
\centering
\includegraphics[scale=0.7]{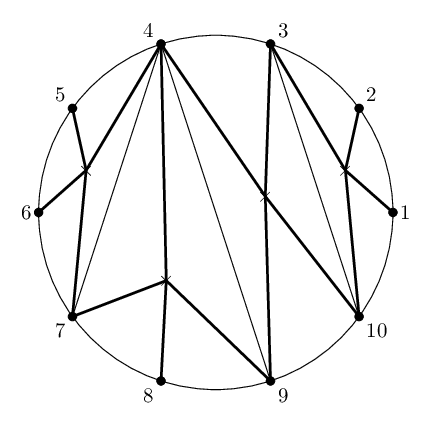}
\caption{An example of an $N$-angulation 
from a saddle-free differential in the case $N=4$ and $n=3$}
\label{saddle-free.fig}
\end{figure}

We write by $\Delta_{\phi}$ the $N$-angulation 
determined by a saddle-free differential $\phi \in \Q$. 
We also write by $h(\delta)$ the horizontal strip which contains a generic trajectory $\delta$. 

Set $\Delta_{\phi}=\{\delta_1,\dots,\delta_n\}$ 
and consider the colored quiver lattice 
\begin{equation*}
L_{\CQ(\Delta_{\phi})}=\oplus_{i=1}^n \Z \alp_i
\end{equation*}
where $\alp_1,\dots,\alp_n$ are the basis 
corresponding to $\delta_1,\dots,\delta_n$. 
Recall from Section \ref{sec_standard_saddle} 
that we can associate standard saddle classes $\gamma_1,\dots,
\gamma_n \in H_{\pm}(\phi)$ to horizontal strips $h(\delta_1),
\dots, h(\delta_ n)$.

\begin{lem}
\label{homology_quiver}
There is an isomorphism of abelian groups
\begin{equation*}
\mu_{\phi} \colon L_{\mathcal{Q}(\Delta_{\phi})}
\to  H_{\pm}(\phi)
\end{equation*}
defined by $\alp_i \mapsto \gamma_i$. The map $\mu_{\phi}$
takes the bilinear form $\left<\,,\,\right>$ on 
$L_{\mathcal{Q}(\Delta_{\phi})}$ to 
the intersection form $I_{\pm}$ on $H_{\pm}(\phi)$.
\end{lem}
{\bf Proof.}
Since standard saddle classes 
$\gamma_1,\dots,\gamma_n$ form the basis of 
$H_{\pm}(\phi)$ by Lemma \ref{standard_saddle}, 
the map $\mu_{\phi}$ is an isomorphism. 

\begin{figure}[!ht]
\centering
\includegraphics[scale=0.7]{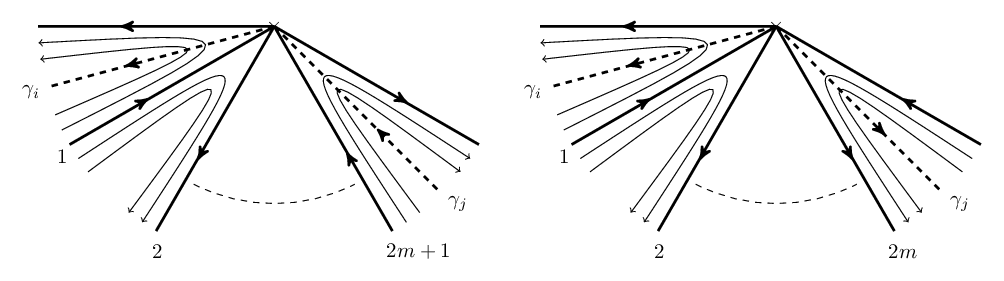}
\caption{Oriented foliations near a zero. }
\label{intersection.fig}
\end{figure}

The remaining part is to prove that 
$\left<\alp_i,\alp_j \right>=I_{\pm}(\gamma_i,\gamma_j)$. 
First note that $\left<\alp_i,\alp_i \right>=I_{\pm}(\gamma_i,\gamma_i)=1+
(-1)^N $ is clear by the definitions of $\left<\,,\,\right>$ and 
$I_{\pm}$. Next consider the case 
$q^{(c)}_{ij} =1$ for some color $c$ and $i \neq j$. 
If $c=2m$, there are $2m+1$ separating trajectories 
counterclockwise from $\gamma_i$ to $\gamma_j$  
as in the left of Figure \ref{intersection.fig}.  
Hence the oriented foliation of the $1$-form 
$\sqrt{\phi}$ is locally given as 
in the left of Figure \ref{intersection.fig}. 
Therefore, 
we have $\left<\alp_i,\alp_j \right>=(-1)^{2m+1}q^{(2m)}_{ij}=-1
=I_{\pm}(\gamma_i,\gamma_j)$. 

Similarly, if $c=2m-1$, the oriented foliation is 
given as in the right of Figure \ref{intersection.fig}, and 
we have $\left<\alp_i,\alp_j \right>=(-1)^{2m}q^{(2m-1)}_{ij}=1
=I_{\pm}(\gamma_i,\gamma_j)$.

\endproof

Let $\phi \in \Q$ be a saddle-free differential and 
$\Delta_{\phi}=\{\delta_1,\dots,\delta_n\}$ 
be the $N$-angulation associated to $\phi$. 
By Corollary \ref{correspondence}, we have 
simple objects $S_1,\dots,S_n \in \H(\Delta_{\phi})$ 
which correspond to $\delta_1,\dots,\delta_n$. 

The next result is the analogue of Lemma 10.5 in \cite{BrSm}. 
\begin{cor}
\label{cor_PC}
There is an isomorphism of 
abelian groups
\begin{equation*}
\nu_{\phi} \colon K(\D^N_n) \cong K(\H(\Delta_{\phi}))\to H_{\pm}(\phi)
\end{equation*}  
defined by $[S_i] \mapsto \gamma_i$. 
The map $\nu_{\phi}$ takes  
the Euler form $\chi$ on $K(\D^N_n)$ to the intersection form $I_{\pm}$ on $H_{\pm}(\phi)$. 
\end{cor}
{\bf Proof.}
First note that by Proposition \ref{correspondence_quiver}, 
we have a canonical isomorphism 
$L_{\CQ(\Delta_{\phi})} \cong L_{\CQ(\H(\Delta_{\phi}))}$. 
Then combining Lemma \ref{quiver_K} and 
Lemma \ref{homology_quiver}, 
the result immediately follows. 
\endproof

\section{From quadratic differentials to stability conditions}
\subsection{Stratifications}
\label{stratification}
Following \cite[Section 5.2]{BrSm}, we introduce the stratification on $\Q$. 
Many results of \cite{BrSm} for the stratification also work  
in our cases. 
For $\phi \in \Q$, denote by $s_{\phi}$ the number of saddle trajectories. 
Define the subsets of $\Q$ by
\begin{equation*}
D_p :=\{\, \phi \in \Q \,\vert\,s_{\phi} \le p\,\}.
\end{equation*}
(Note that our notation differs from that of \cite{BrSm}. Our $D_{p}$ corresponds to 
$B_{2p}$ in \cite{BrSm}.)
Since Lemma \ref{saddle_number} implies that $s_\phi \le n$, 
we have $D_n = \Q$. Note that the stratum $B_0$ 
consists of all saddle-free differentials. 
These subsets have the following property. 
\begin{lem}[\cite{BrSm}, Lemma 5.2]
The subsets $D_p \subset \Q$ form an increasing chain of dense open subsets 
\begin{equation*}
D_0 \subset D_1 \subset \cdots \subset D_n = \Q.
\end{equation*}
\end{lem}

Let $F_0:=D_0$ and $F_p := D_p \,\bs\, D_{p-1}$ for $p \ge 1$. 
Then we have a finite stratification 
\begin{equation*}
\Q = \bigcup_{p=0}^n \,F_p
\end{equation*}
by locally closed subset $F_p$.

\begin{prop}[\cite{BrSm}, Proposition 5.5 and Proposition 5.7]
\label{neighborhood}
Let $p \ge 1$ and assume that $\phi \in F_p$. Then there is a neighborhood $\phi \in U \subset D_p$ and 
a class $0 \neq \alp \in H_{\pm}(\phi)$ such that 
\begin{equation*}
\phi \in U \cap F_p  \Longrightarrow Z_{\phi}(\alp) \in \R.
\end{equation*}
Further, if $p \ge 2$ 
we can take a neighborhood $\phi \in U \subset D_p$ 
which satisfies that $U \cap D_{p-1}$ is connected. 
\end{prop}

Following in \cite[Section 5.2]{BrSm}, 
we define generic differentials. 
$\phi \in \Q$ is called  generic if 
for any $\gamma_1,\gamma_2 \in H_{\pm}(\phi)$, 
\begin{equation*}
\R Z_{\phi}(\gamma_1) = \R Z_{\phi}(\gamma_2) \Longrightarrow \Z \gamma_1 = \Z \gamma_2.
\end{equation*}
It is easy to show that generic differentials are dense in $\Q$.


A differential $\phi \in \Q$ defines the length $l_{\phi}(\gamma)$ 
of a smooth path
$\gamma \colon [0,1] \to \P^1 \,\bs\, \{\infty\}$ by
\begin{equation*}
l_{\phi}(\gamma):=\int_0^1 |\varphi(\gamma(t))|^{1 \slash 2}\left|\frac{d \gamma(t)}{dt} \right| dt
\end{equation*}
where $\varphi(z)$ is given by 
$\phi(z)=\varphi(z) dz^{\otimes 2}$.  
 
Proposition 6.8 in \cite{BrSm} is easily extended to our cases. 
\begin{prop}[\cite{BrSm}, Proposition 6.8]
\label{convergence}
Assume that for a sequence of framed differentials
\begin{equation*}
(\phi_k,\theta_k) \in \Q^{\Gamma}_*, \quad k \ge 1,
\end{equation*} 
the periods $Z_{\phi_k} \circ \theta_k \colon \Gamma \to \C$ converge as $k \to \infty$. 
Further, there is some constant $L>0$, which is independent of $k$, 
such that any saddle connection $\gamma$ of $\phi_k$ satisfies $l_{\phi_k}(\gamma) \ge L$. 
Then there is some subsequence of $\{(\phi_k,\theta_k)\}_{k \ge 1}$ 
which converges in $\Q^{\Gamma}_*$. 
\end{prop}

\subsection{Stability conditions from saddle-free differentials}
Following \cite[Section 10.4]{BrSm}, 
we construct a stability condition on $\D^N_n$ 
from a saddle-free differential $\phi \in \Q$.

\begin{lem}
\label{construct_stability}
$\phi \in B_0 \subset \Q$ be a saddle-free differential. 
Consider the corresponding heart $\H(\Delta_{\phi}) \in \EG_N^{\circ}(\H_0)$ and 
define a linear map 
\begin{equation*}
Z \colon K(\H(\Delta_{\phi})) \lto \C
\end{equation*}
by $Z:=Z_{\phi} \circ \nu_{\phi}$ where $Z_{\phi} \colon H_{\pm}(\phi) \to \C$ is a period of $\phi$ and 
$\nu_{\phi} \colon K(\H(\Delta_{\phi})) \to H_{\pm}(\phi)$ is an isomorphism given by Corollary \ref{cor_PC}.
Then the pair $(Z,\H(\Delta_{\phi}))$ determines 
a stability condition $\sigma(\phi)=(Z,\H(\Delta_{\phi})) \in \Stab^{\circ}(\D_n^N)$. 
\end{lem}
{\bf Proof.}
By Proposition \ref{abel_stability}, 
we only need to check that $Z(E) \in H$ for any 
nonzero object $E \in \H(\Delta_{\phi})$. 
Let $h_1,\dots,h_n$ be horizontal strips determined by $\phi$ 
and $S_1,\dots,S_n$ be corresponding simple objects 
in $\H(\Delta_{\phi})$. 
By the definition of $\nu$ in Corollary \ref{cor_PC}, 
we have $\nu(S_i)=\gamma_i$ where $\gamma_i$ is 
the standard saddle connection associated to the 
horizontal strip $h_i$. As explained in Section \ref{sec_standard_saddle}, 
the period of $\gamma_i$ satisfies 
that $Z_{\phi}(\gamma_i) \in \mathbb{H} \subset H$. 
Hence $Z(S_i) \in H$. 
\endproof

\subsection{Wall crossing}
Let $\phi_1 \in F_1$ be a differential which has just one saddle trajectory. 
Write by $\Delta_{\phi_1}$ the set of $(N-2)$-diagonals obtained by taking one generic trajectory 
from each horizontal strip of $\phi_1$. 
By Lemma \ref{saddle_number}, $n-1$ horizontal strips appear in the decomposition of $\P^1$ by
$\phi_1$, so the number of $(N-2)$-diagonals in $\Delta_{\phi_1}$ is $n-1$.  

By Proposition \ref{neighborhood}, if we take $r >0$ sufficiently small, then 
for any $0 <t \le r$, the differentials 
\begin{equation*}
\phi(t):=e^{+ i \pi t} \cdot \phi_1,\quad \phi^{\sharp}(t):=e^{- i \pi t} \cdot \phi_1
\end{equation*} 
are saddle-free. Write $\phi:=\phi(r)$ and $\phi^{\sharp}:=\phi^{\sharp}(r)$. 

Let $\Delta_{\phi}$ and $\Delta_{\phi^{\sharp}}$ be $N$-angulations 
determined by these saddle-free differentials. 
Write by $\delta_v$ a unique $(N-2)$-diagonal of $\Delta_{\phi}
=\{\delta_1,\dots,\delta_n\}$ 
which is not contained in $\Delta_{\phi_1}$:
\begin{equation*}
\delta_v:=\Delta_{\phi} \,\bs\,\Delta_{\phi_1}.
\end{equation*}

We recall from Section \ref{$N$-angulation} that a new $N$-angulation $\mu^{\sharp}_{\delta_v}(\Delta_{\phi})$ 
is defined by rotating $\delta_v$ one step clockwise.    

\begin{figure}[!ht]
\centering
\includegraphics[scale=0.7]{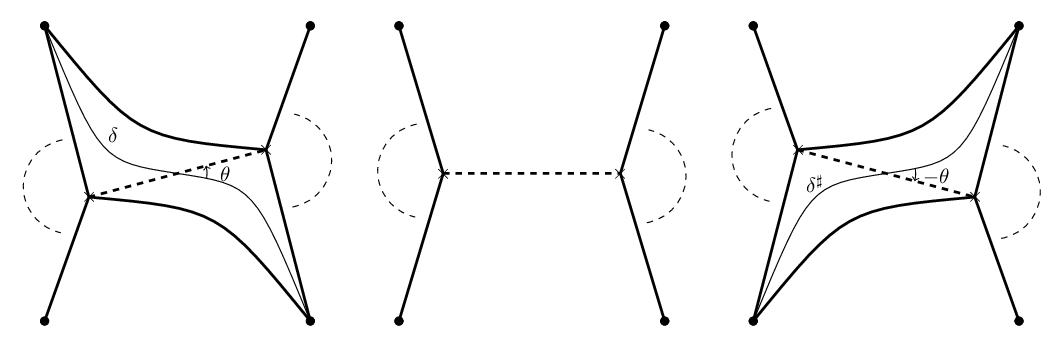}
\caption{Wall crossing.}
\label{WC.fig}
\end{figure}

\begin{lem}
\label{wall_crossing_angulation}
Two $N$-angulations $\Delta_{\phi}$ and $\Delta_{\phi^{\sharp}}$ are related by
\begin{equation*}
\mu^{\sharp}_{\delta_v}(\Delta_{\phi})=\Delta_{\phi^{\sharp}}.
\end{equation*}
\end{lem}
{\bf Proof.} The center of Figure \ref{WC.fig} illustrates a 
unique saddle trajectory of $\phi_1$ and two boundaries of 
half planes or horizontal strips which contain the saddle trajectory. 
 
By Lemma \ref{rotation}, the saddle trajectory maps to 
the saddle connection of phase 
$- \theta$ or $+ \theta$ ($\theta := \d r \slash 2$) for 
$\phi=e^{-i \pi r} \cdot \phi_1$ or 
$\phi^{\sharp}=e^{+i \pi r} \cdot \phi_1$ as 
the left or right of Figure \ref{WC.fig}. 
Comparing the left and right of Figure \ref{WC.fig}, we 
have the result. 

\endproof

By Lemma \ref{ang_mu} and Lemma \ref{wall_crossing_angulation}, 
we have 
\begin{equation*}
\mu^{\sharp}_v \CQ(\Delta_{\phi})= \CQ(\Delta_{\phi^{\sharp}}).
\end{equation*}

As a result, we have a natural isomorphism of colored quiver lattices 
$L_{\mu^{\sharp}_v \CQ(\Delta_{\phi})} \cong L_{\CQ(\Delta_{\phi^{\sharp}})}$. 

For $\phi$ and $\phi^{\sharp}$, we fix the natural basis as follows:
\begin{align*}
L_{\CQ(\Delta_{\phi})} &\cong \bigoplus_{i=1}^n \Z \,\alp_i,   \quad L_{\CQ(\Delta_{\phi^{\sharp}})} \cong \bigoplus_{i=1}^n \Z \,\alp_i^{\prime}, \\
H_{\pm}(\phi) &\cong \bigoplus_{i=1}^n \Z \gamma_i,  \quad \,\,\, H_{\pm}(\phi^{\sharp}) \cong \bigoplus_{i=1}^n \Z \gamma_i^{\prime}
\end{align*}
where $\gamma_1,\dots, \gamma_n$ and $\gamma_1^{\prime},\dots,\gamma_n^{\prime}$ are standard saddle classes. 

Let 
\begin{equation*}
F_v \colon L_{\CQ(\Delta_{\phi^{\sharp}})} \lto L_{\CQ(\Delta)}
\end{equation*}
be the linear isomorphism defined in Section \ref{sec_mutation_K}. 

\begin{figure}[!ht]
\centering
\includegraphics[scale=0.8]{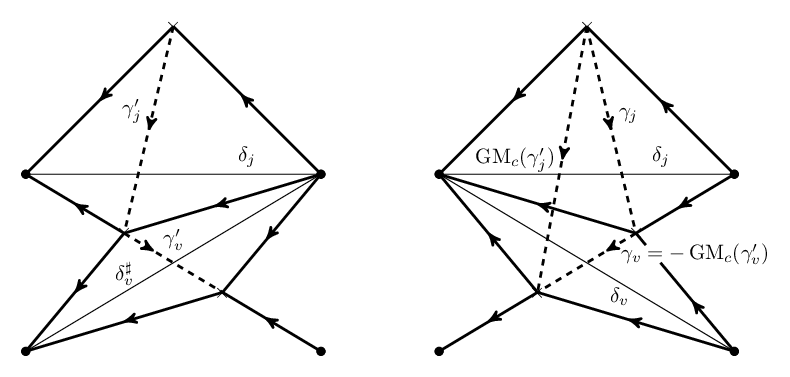}
\caption{Changes of cycles along the path $c$. }
\label{GM.fig}
\end{figure}

The relationship between $F_v$ and the Gauss-Manin connection 
\begin{equation*}
\GM_c \colon H_{\pm}(\phi^{\sharp}) \lto H_{\pm}(\phi)
\end{equation*}
along the path $c(t)=e^{i \pi t} \cdot \phi_1\,(t \in [-r,r])$ is given by the 
following diagram. 
\begin{lem}
\label{diagram2}
Let $\mu_{\phi} \colon L_{\CQ(\Delta_{\phi})} \to H_{\pm}(\phi)$ and 
$\mu_{\phi^{\sharp}} \colon L_{\CQ(\Delta_{\phi^{\sharp}})} \to H_{\pm}(\phi^{\sharp})$ be the isomorphisms given by 
Lemma \ref{homology_quiver}. 
Then the diagram
\begin{equation*}
\xymatrix{ 
L_{\CQ(\Delta_{\phi^{\sharp}})} \ar[d]_{\mu_{\phi^{\sharp}}} \ar[rr]^{F_v} && 
L_{\CQ(\Delta_{\phi})} \ar[d]^{\mu_{\phi}} \\ 
H_{\pm}(\phi^{\sharp}) \ar[rr]^{\GM_c} && H_{\pm}(\phi)
}
\end{equation*}
commutes. 
\end{lem}
{\bf Proof.} As in Figure \ref{GM.fig} 
if $q^{(0)}_{vj}=1$, then we have 
\begin{equation*}
\GM_c(\gamma_j^{\prime})=\gamma_j + \gamma_v, \quad \GM_c(\gamma_v^{\prime})=-\gamma_v.
\end{equation*}
On the other hand if $q^{(0)}_{jv}=0$, then we can check $\GM_c(\gamma_j^{\prime})=\gamma_j$. 
Since $\mu_{\phi}(\alp_j)=\gamma_j$ and $\mu_{\phi^{\sharp}}(\alp_j^{\prime})=\gamma_j^{\prime}$, 
the result follows.

\endproof

For $0 <t \le r$, we denote by
\begin{equation*}
\sigma(t):=\sigma(\phi(t)), \quad \sigma^{\sharp}(t):=\sigma(\phi^{\sharp}(t))
\end{equation*}
stability conditions in $\Stab^{\circ}(\D_n^N)$ constructed by Lemma \ref{construct_stability} 
for saddle-free differentials $\phi(t)$ and $\phi^{\sharp}(t)$. 

The following result is the analogue of Proposition 10.7 in \cite{BrSm}.

\begin{prop}
\label{diagram3}
Let $\nu_{\phi} \colon K(\D_n^N) \to H_{\pm}(\phi)$ and $\nu_{\phi^{\sharp}} \colon K(\D_n^N) \to H_{\pm}(\phi^{\sharp})$ 
be the isomorphisms given by Corollary \ref{cor_PC}. 
Then there is a unique spherical twist $\Phi \in \Sph(\D_n^N)$ such that the following $(1)$ and $(2)$ hold.
\begin{itemize}
\item[(1)] The diagram 
\begin{equation*}
\xymatrix{ 
K(\D_n^N) \ar[d]_{\nu_{\phi^{\sharp}}} \ar[rr]^{[\Phi]} && 
K(\D_n^N) \ar[d]^{\nu_{\phi}} \\ 
H_{\pm}(\phi^{\sharp}) \ar[rr]^{\GM_c} && H_{\pm}(\phi)
}
\end{equation*}
commutes. 
\item[(2)] 
\begin{equation*}
\lim_{t \to +0} \sigma(t)= \lim_{t \to +0}\Phi \cdot 
\sigma^{\sharp}(t).
\end{equation*}
\end{itemize}
\end{prop}
{\bf Proof.} Part $(1)$ immediately follows from Lemma \ref{diagram11} and Lemma \ref{diagram2}. 
For the proof of part $(2)$, we first note that 
by Lemma \ref{lem_spherical} and Lemma \ref{wall_crossing_angulation}, there is a unique spherical twist $\Phi \in \Sph(\D_n^N)$ such that 
\begin{equation*}
\Phi \cdot \H(\Delta_{\phi^{\sharp}})=
\H(\Delta_{\phi})_{S_v}^{\sharp}.
\end{equation*}

By the definition of $\sigma(t)$ and $\sigma^{\sharp}(t)$ in Lemma \ref{construct_stability}, 
$\sigma(t)$ is contained in $\Stab(\H(\Delta_{\phi}))$ and 
$\Phi \cdot \sigma^{\sharp}(t)$ is contained in $\Stab(\H(\Delta_{\phi})_{S_v}^{\sharp})$. 
Further by Part $(1)$, central charges of 
$\sigma(t)$ and $\Phi \cdot \sigma^{\sharp}(t)$ approach 
the same point in $\Hom(K(\D_n^N),\C)$ as $t \to +0$.  
Hence by Lemma \ref{gluing}, we can take some open subset 
\begin{equation*}
U \subset \Stab(\H(\Delta_{\phi})) 
\cup \Stab(\H(\Delta_{\phi})_{S_v}^{\sharp})
\end{equation*}
such that $U$ is mapped isomorphically on 
$\Hom(K(\D_n^N),\C)$, and $U$ contains both 
$\sigma(t)$ and $\Phi \cdot \sigma^{\sharp}(t)$. 
This implies the result. 
\endproof

\subsection{Construction of isomorphisms}
\begin{defi}
\label{kappa}
Let $(\phi,\theta) \in \Q^{\Gamma}$ be 
a framed saddle-free differential and 
$\nu_{\phi} \colon K(\D_n^N) \to H_{\pm}(\phi)$ be 
an isomorphism given by Corollary \ref{cor_PC}. 
We define the isomorphism $\kappa_{\phi}$ by the composition of 
$\theta$ and $\nu_{\phi}^{-1}$: 
\begin{equation*}
\kappa_{\phi}:=\nu_{\phi}^{-1} \circ \theta \colon 
\Gamma \to K(\H(\Delta_{\phi})) \cong K(\D_n^N).
\end{equation*}
\end{defi}

Let $(\phi_0,\theta_0) \in \Q^{\Gamma}$ be the 
fixed framed differential in Section \ref{sec_local_system}
and $\Q_*^{\Gamma}$ be the connected 
component which contains $(\phi_0,\theta_0)$. 
In the following, we also assume that 
$(\phi_0,\theta_0)$ is saddle-free. 

We identify $K(\D_n^N)$ with 
$\Gamma$ by using the isomorphism 
$\kappa_0:=\kappa_{\phi_0}$ determined 
by Definition \ref{kappa}:
\begin{equation*}
\kappa_0 \colon \Gamma \iso K(\D_n^N). 
\end{equation*}

\begin{lem}
\label{lem_kappa}
Let $(\phi,\theta) \in \Q_*^{\Gamma}$ be 
a framed saddle-free differential and 
$\kappa_{\phi} \colon \Gamma \to K(\D_n^N)$ be an isomorphism 
given by Definition \ref{kappa}. 
Then there is a spherical twist $\Phi \in \Sph(\D_n^N)$ 
such that the diagram 
\begin{equation*}
\xymatrix{ 
& \Gamma \ar[dl]_{\kappa_{\phi}} \ar[dr]^{\kappa_0}  & \\ 
K(\D_n^N) \ar[rr]^{[\Phi]} && K(\D_n^N).
}
\end{equation*}
commutes. In addition, if there are two such spherical twists 
$\Phi_1$ and $\Phi_2$, then 
the induced linear isomorphisms on $K(\D_n^N)$ are the same:
\begin{equation*}
[\Phi_1]=[\Phi_2].
\end{equation*}
\end{lem}
{\bf Proof.} $[\Phi_1]=[\Phi_2]$ immediately follows from 
\begin{equation*}
[\Phi_1] \circ \kappa_{\phi} = \kappa_0 = [\Phi_2] 
\circ \kappa_{\phi}.
\end{equation*}
Therefore, we only need to show that such a spherical twist exists. 
Take a path 
$c \colon [0,1] \to \Q_*^{\Gamma}$ from $(\phi,\theta)$ to 
$(\phi_0,\theta_0)$. 
By Proposition \ref{neighborhood}, we can deform the path $c$ 
to lie in $D_1$ and to intersect only finitely many points 
of $F_1$, which consists of framed differentials with just one 
saddle trajectory. 

By applying Proposition \ref{diagram3} to each of these points,  
we have 
a spherical twist $\Phi \in \Sph(\D_n^N)$ fitting into a diagram
\begin{equation*}
\xymatrix{ 
K(\D_n^N) \ar[d]_{\nu_{\phi}} \ar[rr]^{[\Phi]} && 
K(\D_n^N) \ar[d]^{\nu_0} \\ 
H_{\pm}(\phi) \ar[rr]^{\GM_c} && H_{\pm}(\phi_0)
}
\end{equation*}
where $\nu_{\phi}$ and $\nu_0:=\nu_{\phi_0}$ are 
isomorphisms given by Corollary \ref{cor_PC}. 
Further, recall the diagram 
\begin{equation*}
\xymatrix{ 
& \Gamma \ar[dl]_{\theta} \ar[dr]^{\theta_0}  & \\ 
H_{\pm}(\phi) \ar[rr]^{\GM_c} && H_{\pm}(\phi_0)
}
\end{equation*}
in Section \ref{sec_local_system}. Combining the above two diagrams, 
we have the result. 
\endproof

Let $\Sph^0(\D^N_n) \subset \Sph(\D^N_n)$ be the subgroup 
of spherical twists which act as the identity on $K(\D^N_n)$. 
We recall that $\Sph(\D_n^N) \cong B_{n+1}$ 
by Theorem \ref{faithful}, 
and $\Sph(\D_n^N)$ acts as the group $W_{\pm}$ on 
$K(\D^N_n)$ (see Section \ref{sec_spherical}).  
Hence the group $\Sph^0(\D^N_n)$ is 
isomorphic to the group $P_{\pm}$ in Section \ref{sec_braid}.

By completely the same argument of Proposition 11.3 in \cite{BrSm}, we can construct 
the following equivariant holomorphic map.  
\begin{prop}
\label{map}
There is a $W_{\pm}$-equivariant and $\C$-equivariant 
holomorphic map of complex manifolds $K_N$ such that 
the diagram
\begin{equation*}
\xymatrix{
\Q_*^{\Gamma} \ar[rr]^{K_N} \ar[dr]_{\mathcal{W}_N} && 
\Stab^{\circ}(\D^N_n) \slash  \Sph^0(\D_n^N) \ar[dl]^{\mathcal{Z}}  \\
& \Hom_{\Z}(\Gamma,\C) & 
}
\end{equation*}
commutes where  $\mathcal{W}_N$ and $\mathcal{Z}$ are 
local isomorphisms 
in Proposition \ref{local_iso_period} 
and Theorem \ref{localiso} respectively.
\end{prop}
{\bf Proof.} First we see the $\C$-equivariant property of $K_N$. 
Define the action of $t \in \C$ on $\Hom_{\Z}(\Gamma,\C)$ 
by $Z \mapsto e^{-i \pi t} \cdot Z$ 
for $Z \in \Hom_{\Z}(\Gamma,\C)$.  
Then by Definition \ref{C_action} and Lemma \ref{rotation}, 
the map $\mathcal{W}_N$ is $\C$-equivariant. 
On the other hand, by the definition of $\C$-action on 
$\Stab^{\circ}(\D^N_n) $ given in Section \ref{sec_group},  
the map $\mathcal{Z}$ is also $\C$-equivariant. 
Since both $\mathcal{W}_N$ and $\mathcal{Z}$ are 
local isomorphisms by Proposition \ref{local_iso_period} and Theorem \ref{localiso}, 
if $K_N$ exists then $K_N$ is $\C$-equivariant. 

Next we construct the map $K_N$. 
Recall from Section \ref{stratification} that there is 
a stratification 
\begin{equation*}
D_0 \subset D_1 \subset \cdots \subset D_n = \Q.
\end{equation*}
We extend the stratification of $\Q$ to $\Q_*^{\Gamma}$ by 
the obvious way (and use the same notation). 

Here we give the construction of $K_N$ only on the stratum 
$D_0$ which consists of framed saddle-free differentials. 
By using Proposition \ref{neighborhood} and 
Proposition \ref{diagram3}, 
we can extend $K_N$ on 
the stratum $D_p \,(p \ge 1)$ by doing  
the same argument of Proposition 11.3 in \cite{BrSm}.

Let $(\phi,\theta) \in D_0$ be a framed saddle-free differential. 
By Lemma \ref{construct_stability}, we can construct 
a stability condition $\sigma(\phi) \in \Stab^{\circ}(\D_n^N)$. 
Let $\Phi \in \Sph(\D^N_n)$ be a spherical twist 
 for $(\phi,\theta)$ determined by Lemma \ref{lem_kappa}. 
We note that $\Phi$ is determined up to the action of 
$\Sph^0(\D^N_n) \subset \Sph(\D_n^N)$. 
Therefore, if we define 
the map $K_N \colon D_0 \to \Stab^{\circ}(\D_n^N) 
\slash \Sph^0(\D^N_n)$ 
by
\begin{equation*}
K_N((\phi,\theta)):=\Phi ( \sigma(\phi)) \in 
\Stab^{\circ}(\D_n^N) \slash \Sph^0(\D^N_n),
\end{equation*}  
then $K_N$ is well-defined.

Finally we show that the map $K_N$ is $W_{\pm}$-equivariant. 
Note that the action of $W_{\pm}$ preserves the 
stratification of $\Q_*^{\Gamma}$ 
and commutes with the action of $\C$. 
Hence it is sufficient to prove that $K_N$ is 
$W_{\pm}$-equivariant only on $D_0$,  
and this is clear by the construction of $K_N$ on $D_0$. 
\endproof

Let $(\phi,\theta) \in \Q_*^{\Gamma}$ and
we set $\sigma =K_N((\phi,\theta))$. 
The stability condition $\sigma$ is determined up to 
the action of $\Sph^0(\D^N_n) \subset \Aut(\D^N_n)$. 
For a class $0 \neq \alp \in \Gamma$, 
let $\mathcal{M}_\sigma(\alp)$ be the moduli space of 
$\sigma$-stable objects which have the class $\alp$ 
and the phase in $(0,1]$ 
(here, we treat $\mathcal{M}_\sigma(\alp)$ 
as a set).   


The same argument of Proposition 11.11 in \cite{BrSm} 
gives the following. 
\begin{prop}
\label{iso}
The map $K_N$ in Proposition \ref{map} is an isomorphism. 
\end{prop}
{\bf Proof.} Here we recall 
the argument of Proposition 11.11 in \cite{BrSm}. 

We first show that $K_N$ is injective. 
Let $(\phi_1,\theta_1), (\phi_2,\theta_2) \in \Q^{\Gamma}_*$ 
such that $K_N((\phi_1,\theta_1))=K_N((\phi_2,\theta_2))$. 
Since the maps $\mathcal{W}_N$ and $\mathcal{Z}$ 
are local isomorphisms  
by Proposition \ref{local_iso_period} and 
Theorem \ref{localiso} and $K_N$ commutes with 
these two maps, 
we can deform 
$\phi_1$ and $\phi_2$ to be saddle-free differentials. 

By the construction of $K_N$ in Proposition \ref{map}, 
we can write $K_N((\phi_i,\theta_i))= \Phi_i( \sigma (\phi_i))$ 
for some $\Phi_i \in \Sph(\D_n^N)$ 
determined by Lemma \ref{lem_kappa}. 
Since $ \Phi_1( \sigma (\phi_1))$ and 
$ \Phi_2( \sigma (\phi_2))$ have the same heart 
(up to the action of $\Sph^0(\D_n^N)$), 
we have 
\begin{equation*}
\Phi_1 \cdot \H(\Delta_{\phi_1})  \equiv
\Phi_2 \cdot \H(\Delta_{\phi_2}) \mod \Sph^0(\D^N_n).
\end{equation*}
Recall Remark \ref{heart_fundamental} that $\H(\Delta_{\phi_1})$ and $\H(\Delta_{\phi_2})$ 
are in the fundamental domain 
$\EG_N^{\circ}(\H_0)$. Hence Theorem \ref{fundamental} 
implies that $\H(\Delta_{\phi_1})=\H(\Delta_{\phi_2})$ 
and $[\Phi_1]=[\Phi_2]$ as linear maps on $K(\D_n^N)$. 
As a result, we have $\Delta_{\phi_1}=\Delta_{\phi_2}$ by 
Corollary \ref{correspondence}. 
Thus $\phi_1$ and $\phi_2$ have the same $N$-angulations 
and the same periods. Then Proposition 4.9 in \cite{BrSm} 
implies that $\phi_1=\phi_2$. 
On the other hand, $[\Phi_1]=[\Phi_2]$ together with the diagram of Lemma \ref{lem_kappa} implies  
that $\theta_1=\theta_2$. 

Next to prove the surjectivity of $K_N$, we show that 
the image of $K_N$ is open and closed. 
The injectivity of $K_N$ together with the local isomorphism property 
of $\mathcal{W}_N$ and $\mathcal{Z}$ implies that 
the image of $K_N$ is open. 

Set $\sigma_k:=K_N((\phi_k,\theta_k)) 
\in \Stab^{\circ}(\D_n^N) \slash \Sph^0(\D^N_n)$ 
for the sequence of 
framed differentials $\{(\phi_k,\theta_k)\}_{k \ge 1}$ 
in $\Q_*^{\Gamma}$, 
and assume that $\sigma_k$ converge to the point $\sigma_{\infty}$ 
as $k \to \infty$. 
We can assume that $(\phi_k,\theta_k)$ 
is generic saddle-free since generic 
saddle-free differentials are dense in $\Q^{\Gamma}_*$. 
In this setting, we show that there is a subsequence of 
$\{(\phi_k,\theta_k)\}_{k \ge 1}$ which converges 
some point in $ \Q_*^{\Gamma}$. 
First note that the periods $Z_{\phi_k} \circ \theta_k$ converge 
as $k \to \infty$ 
since the central charge $Z_k$ of $\sigma_k$ also converge as 
$k \to \infty$ by the assumption.  
Set 
\begin{equation*}
M_k:=\inf \{|Z_k(E)| \in \R_{\ge 0} \,\vert\, [E] \in 
\class^{\mathrm{ss}}(\sigma_k)\} \quad\text{for }k=1,\dots,\infty
\end{equation*}
where $Z_k$ is the central charge of $\sigma_k$. 
The support property of $\sigma_k$ implies 
that $M_k >0$. Since $M_k \to M_{\infty}$ as $k \to \infty$, 
there is some constant $L >0$ such that 
$M_k \ge L$ for all $k$. 
By Theorem 11.6 and Proposition 11.7 in \cite{BrSm}, 
the length $l_{\phi_k}(\gamma)$ of the saddle connection 
$\gamma$ for $\phi_k$ coincides with the mass 
$|Z_k(S_{\gamma})|$ of the corresponding stable object $S_{\gamma}$. 
Thus we can apply 
Proposition \ref{convergence} for 
$\{(\phi_k,\theta_k)\}_{k \ge 1}$, and obtain 
the convergent subsequence.   
\endproof

\begin{thm}
\label{thm_main}
There is a $B_{n+1}$-equivariant and 
$\C$-equivariant isomorphism of complex manifolds 
$\widetilde{K_N}$ such that the diagram
\begin{equation*}
\xymatrix{
\widetilde{M_n} \ar[rr]^{\widetilde{K_N}} \ar[dr]_{\mathcal{W}_N} && 
\Stab^{\circ}(\D^N_n)  \ar[dl]^{\mathcal{Z}}  \\
& \Hom_{\Z}(\Gamma,\C) &
}
\end{equation*}
commutes. 
\end{thm}
{\bf Proof.} First recall from Corollary \ref{universal} that 
$\widetilde{M_n}$ is the universal cover of the space 
$\Q_*^{\Gamma}$ and 
$\widetilde{M_n} \slash P_{\pm} \iso \Q_*^{\Gamma}$. 
Further,  by Proposition \ref{map} and 
Proposition \ref{iso}, we have the isomorphism 
\begin{equation*}
\widetilde{M_n} \slash P_{\pm} \iso \Stab^{\circ}(\D^N_n) 
\slash \Sph^0(\D_n^N).
\end{equation*}  
Since $\Sph^0(\D_n^N) \cong P_{\pm}$ and the quotient map 
$\Stab^{\circ}(\D^N_n) \to \Stab^{\circ}(\D^N_n) \slash 
 \Sph^0(\D_n^N)$ is a covering map by Lemma \ref{free},  
we can lift the map 
$K_N$ to the $B_{n+1}$-equivariant isomorphism 
\begin{equation*} 
\widetilde{K_N} \colon 
\widetilde{M_n} \iso \Stab^{\circ}(\D_n^N). 
\end{equation*}
\endproof

\begin{cor}
\label{contractible}
The distinguished connected component 
$\Stab^{\circ}(\D_n^N)$ is contractible. 
\end{cor}
{\bf Proof.}
Since $\widetilde{M_n}$ is contractible 
(for example, see Lemma 1.28 in \cite{KT}), the result follows. 
\endproof 

\subsection{Autoequivalence groups}
In this section, we determine the group of 
autoequivalences $\Aut^{\circ}(\D_n^N)$ which 
preserves the distinguished connected component of 
the exchange graph $\EG^{\circ}(\D_n^N)$. 

\begin{lem}
\label{nil}
Let $\Phi \in \Aut^{\circ}(\D_n^N)$ an autoequivalence which 
acts trivially 
on $\EG^{\circ}(\D_n^N)$. 
Then $\Phi \cong \id$.
\end{lem}
{\bf Proof. } 
Since $\Phi$ preserves all labeled edges in 
$\EG^{\circ}(\D_n^N)$, it preserves simple objects $R_1,\dots,
R_n$ in the standard heart $\H_0$, i.e. $R_i \cong \Phi(R_i)$. 
We note that the dimension of the homomorphism space between $R_i$ and $R_j$ is at most one because the number of 
arrows between 
two vertices in the $A_n$-quiver is at most one. 
Hence $\Phi$ preserves 
the homomorphism space $\Hom_{D_n^N}(R_i,R_j)$. Since 
$\D_n^N$ is generated by $R_1,\dots,R_n$ and their shifts,  
$\Phi$ is isomorphic to the identity functor.
\endproof 

Let $\Pi_{\d}$ be a $\d$-gon as in Section \ref{$N$-angulation}. 
The mapping class group of $\Pi_{\d}$ is given by a  
cyclic group $\Z \slash \d \Z$. We assume that the generator 
$1 \in \Z \slash \d \Z$ corresponds to one step clockwise rotation of $\Pi_{\d}$.

\begin{thm}[\cite{BrSm}, Theorem 9.10]
\label{exact}
Assume that $n \ge 2$. Then 
there is a short exact sequence of groups
\begin{equation*}
1 \to \Sph(\D_n^N) \to \Aut^{\circ}(\D_n^N) 
\to \Z \slash \d \,\Z \to 1.
\end{equation*}
\end{thm}
{\bf Proof. } The same argument of Theorem 9.10 in \cite{BrSm} 
gives the following short exact sequence:
\begin{equation*}
1 \to \Sph(\D_n^N) \to \Aut^{\circ}(\D_n^N) \slash  
\mathrm{Nil}^{\circ}(\D_n^N)
\to \Z \slash \d \,\Z \to 1
\end{equation*}
where $\mathrm{Nil}^{\circ}(\D_n^N)$ is the subgroup  
of  $\Aut^{\circ}(\D_n^N)$ consisting 
of autoequivalences which act trivially on $\EG^{\circ}(\D_n^N)$. 
Lemma \ref{nil} implies that $\mathrm{Nil}^{\circ}(\D_n^N)=1 $.
\endproof

For an integer $m \neq 0$, let $\Z[m] \subset \Aut^{\circ}(\D_n^N)$ be the infinite cyclic 
group generated by the degree $m$ shift functor 
$[m] \in \Aut^{\circ}(\D_n^N)$. 

Recall from Theorem \ref{faithful} that 
there is an isomorphism 
\begin{equation*}
B_{n+1} \iso \Sph(\D_n^N), \quad 
\sigma_i \mapsto \Phi_{S_i}. 
\end{equation*}
By Lemma 4.14 in \cite{ST}, 
the generator of center 
$C \in B_{n+1}$ (see Section \ref{sec_braid}) is 
mapped to the shift functor $[-\d] \in \Sph(\D_n^N)$ under 
the above isomorphism.

Define the diagonal action of the element 
$\d \in \d\,\Z$ on $B_{n+1} \times \Z$ by 
\begin{equation*}
(\sigma, l ) \mapsto (\sigma \cdot C, l+\d)
\quad \text{for }(\sigma,l) \in B_{n+1} \times \Z. 
\end{equation*}

\begin{thm}
Assume that $n \ge 2$. 
Then the group $\Aut^{\circ}(\D^N_n)$ is given by
\begin{equation*}
\Aut^{\circ}(\D^N_n) \cong (B_{n+1} \times \Z) \,\slash \,\d\,\Z
\end{equation*}
where $1 \in \Z$ acts as the shift functor 
$[1] \in \Aut^{\circ}(\D_n^N)$.
\end{thm}
{\bf Proof. } Consider the group homomorphism 
\begin{equation*}
\rho \colon B_{n+1} \times \Z \lto \Aut^{\circ}(\D_n^N)
\end{equation*}
defined by $(\sigma_i,l) \mapsto \Phi_{S_i}[l]$. 

By Definition \ref{C_action}, $1 \in \C$ acts as a rotation of 
the phase $2 \pi \slash \d$ in the clockwise direction 
on $\Q$. This implies 
that single shift functor $[1] \in \Aut^{\circ}(\D_n^N)$ 
is mapped to the generator 
$1 \in \Z \slash \d \Z$ under the map of the 
short exact sequence in Theorem \ref{exact}.  
Hence $\rho$ is surjective. 
On the other hand, the kernel of $\rho$ is given by 
\begin{equation*}
\{\,(C^m, m \d) \in B_{n+1} \times \Z\,\vert\,m \in \Z\} 
\end{equation*}
since $\rho(C)=[-\d]$. 
\endproof

\begin{rem}
In the case $n=1$, the direct calculation gives that  
$\Sph(\D_1^N) \cong \Z[-N+1]$ and 
$\Aut^{\circ}(\D_1^N) \cong \Z[1]$. 

\end{rem}


\begin{thebibliography}{BBD82}

\bibitem[Ami09]{Am}
C~Amiot.
\newblock Cluster categories for algebras of global dimension 2 and quivers
  with potential.
\newblock {\em Ann. Inst. Fourier (Grenoble)}, 59(6):2525--2590, 2009.

\bibitem[BBD82]{BBD}
A.~A. Be{\u\i}linson, J.~Bernstein, and P.~Deligne.
\newblock Faisceaux pervers.
\newblock In {\em Analysis and topology on singular spaces, {I} ({L}uminy,
  1981)}, volume 100 of {\em Ast\'erisque}, pages 5--171. Soc. Math. France,
  Paris, 1982.

\bibitem[BM08]{BaMa}
K.~Baur and R.~J. Marsh.
\newblock A geometric description of {$m$}-cluster categories.
\newblock {\em Trans. Amer. Math. Soc.}, 360(11):5789--5803, 2008.

\bibitem[Bri06]{Br5}
T.~Bridgeland.
\newblock Stability conditions on a non-compact {C}alabi-{Y}au threefold.
\newblock {\em Comm. Math. Phys.}, 266(3):715--733, 2006.

\bibitem[Bri07]{Br1}
T.~Bridgeland.
\newblock Stability conditions on triangulated categories.
\newblock {\em Ann. of Math.}, 166(2):317--345, 2007.

\bibitem[Bri09]{Br4}
T.~Bridgeland.
\newblock Spaces of stability conditions.
\newblock In {\em Algebraic geometry---{S}eattle 2005. {P}art 1}, volume~80 of
  {\em Proc. Sympos. Pure Math.}, pages 1--21. Amer. Math. Soc., Providence,
  RI, 2009.



\bibitem[BS15]{BrSm}
T.~Bridgeland and I.~Smith.
\newblock Quadratic differentials as stability conditions.
\newblock {\em Publ. Math. de l'IHES.}, 121:155--278, 2015.


\bibitem[BQS]{BQS}
T.~Bridgeland, Y.~Qiu and T.~Sutherland.
\newblock Stability conditions and the $A_2$ quiver.
\newblock arXiv:1406.2566.


\bibitem[BT09]{BT}
A.~B. Buan and H.~Thomas.
\newblock Coloured quiver mutation for higher cluster categories.
\newblock {\em Adv. Math.}, 222(3):971--995, 2009.

\bibitem[Dou02]{Dou}
R.~Douglas.
\newblock Dirichlet branes, homological mirror symmetry, and stability.
\newblock In {\em Proceedings of the {I}nternational {C}ongress of
  {M}athematicians, {V}ol. {III} ({B}eijing, 2002)}, pages 395--408, Beijing,
  2002. Higher Ed. Press.

\bibitem[Dub04]{Dub1}
B.~Dubrovin.
\newblock On almost duality for {F}robenius manifolds.
\newblock In {\em Geometry, topology, and mathematical physics}, volume 212 of  {\em Amer. Math. Soc. Transl. Ser. 2.}, pages 75--132. Amer. Math. Soc, Providence, RI, 2004.

\bibitem[Gin]{Ginz}
V.~G. Ginzburg.
\newblock Calabi-{Y}au algebras.
\newblock arXiv:math/0612139.

\bibitem[GMN]{GMN}
D.~Gaiotto, G.~Moore, and A.~Neitzke.
\newblock Wall-crossing, {H}itchin {S}ystems, and the {WKB} {A}pproximation.
\newblock arXiv:0907.3987.

\bibitem[HRS96]{HRS}
D.~Happel, I.~Reiten, and S.~O. Smal{\o}.
\newblock Tilting in abelian categories and quasitilted algebras.
\newblock {\em Mem. Amer. Math. Soc.}, 120(575):viii+ 88, 1996.

\bibitem[IUU10]{IUU}
A.~Ishii, K.~Ueda, and H.~Uehara.
\newblock Stability conditions on {$A_n$}-singularities.
\newblock {\em J. Differential Geom.}, 84(1):87--126, 2010.

\bibitem[Kel11]{Kel2}
B~Keller.
\newblock Deformed {C}alabi-{Y}au completions.
\newblock {\em J. Reine Angew. Math.}, 654:125--180, 2011.
\newblock With an appendix by Michel Van den Bergh.

\bibitem[KN13]{KN}
B.~Keller and P.~Nicol\'as.
\newblock Weight structures and simple dg modules for positive dg algebras
\newblock {\em Int. Math. Res. Not.}, 1028--1078, 2013.


\bibitem[KQ]{KQ}
K.~King and Y.~Qiu.
\newblock Exchange graphs and {E}xt quivers.
\newblock to appear in {\em Adv. Math.}, arXiv:1109.2924v3.  
 
 
 
\bibitem[KS]{KS}
M.~Kontsevich and Y.~Soibelman.
\newblock Stability structures, motivic {D}onaldson-{T}homas invariants and
  cluster transformations.
\newblock arXiv:0811.2435.

\bibitem[KS02]{KhSe}
M.~Khovanov and P.~Seidel.
\newblock Quivers, floer cohomology, and braid group actions.
\newblock {\em J. Amer. Math. Soc.}, 15:203--271, 2002.

\bibitem[KT08]{KT}
C.~Kassel and V.~Turaev.
\newblock {\em Braid groups}, volume 247 of {\em Graduate Texts in
  Mathematics}.
\newblock Springer, New York, 2008.
\newblock With the graphical assistance of Olivier Dodane.

\bibitem[KY11]{KY}
B.~Keller and D.~Yang.
\newblock Derived equivalences from mutations of quivers with potential.
\newblock {\em Adv. Math.}, 226(3):2118--2168, 2011.

\bibitem[Qiu15]{Q}
Y.~Qiu.
\newblock Stability conditions and quantum dilogarithm identities for {D}ynkin
  quivers.
\newblock {\em Adv. Math.}, 269:220--264, 2015.


\bibitem[Sai83]{Sai}
K.~Saito.
\newblock Period Mapping Associated to a Primitive Form. 
\newblock In {\em Publ. RIMS.} 19:1231--1264, 1983.

\bibitem[ST01]{ST}
P.~Seidel and R.~P. Thomas.
\newblock Braid group actions on derived categories of coherent sheaves.
\newblock {\em Duke Math. J.}, 108(1):37--108, 2001.

\bibitem[Str84]{Str}
K~Strebel.
\newblock {\em Quadratic differentials}, volume~5 of {\em Ergebnisse der
  Mathematik und ihrer Grenzgebiete (3) [Results in Mathematics and Related
  Areas (3)]}.
\newblock Springer-Verlag, Berlin, 1984.

\bibitem[Sut]{Suther}
T.~Sutherland.
\newblock The modular curve as the space of stability conditions of a {C}{Y}3
  algebra.
\newblock arXiv:1111.4184.

\bibitem[Tho06]{T}
R.~P. Thomas.
\newblock Stability conditions and the braid group.
\newblock {\em Comm. Anal. Geom.}, 14(1):135--161, 2006.

\bibitem[TY02]{TY}
R.P. Thomas and S.T. Yau.
\newblock Special {L}agrangians, stable bundles and mean curvature flow.
\newblock {\em Comm. Anal. Geom.}, 10:1075--1113, 2002.

\bibitem[Woo10]{Wol}
J.~Woolf.
\newblock Stability conditions, torsion theories and tilting.
\newblock {\em J. Lon. Math. Soc.}, 82:663--682, 2010.

\end{thebibliography}
\end{document}